\documentclass[fleqn,oneside,a4paper,12pt]{article}
\title{Friezes of type \tilde{D}}
\usepackage[english]{babel}
\usepackage[T1]{fontenc}
\usepackage[a4paper,tmargin=3cm,bmargin=3cm,rmargin=2.1cm,lmargin=2.1cm]{geometry}
\usepackage{amsfonts}
\usepackage{qsymbols}
\usepackage{amsmath}
\usepackage{euscript}
\usepackage{wasysym}
\usepackage{babel}
\usepackage{amssymb}
\usepackage[T1]{fontenc}
\usepackage[all]{xy}
\usepackage{graphicx}
\usepackage{epsf}
\usepackage{phonetic}
\usepackage{url}
\usepackage{showkeys}
\usepackage{stmaryrd}
\usepackage{mathrsfs}
\usepackage{expdlist}
\usepackage{fancyhdr}
\usepackage{theorem,amsfonts}
\theoremstyle{break}

\newtheorem{defi}{Definition}[section]
\newtheorem{exem}[defi]{Example}

\newtheorem{rem}[defi]{Remark}
\newtheorem{prop}[defi]{\textsc{Proposition}}
\newtheorem{lem}[defi]{Lemma}
\newtheorem{theo}[defi]{Theorem}
\newtheorem{Cor}[defi]{Corollary}

\title{Friezes of type $\tilde{\mathbb{D}}$}
\author{Kodjo Essonana Magnani}
\begin{document}
\maketitle
\begin{abstract}
In this article, we establish a relation between the values of a frieze of type $\tilde{\mathbb{D}}_{n}$ and some values of an $SL_{2}$-tiling $t$ associated with a particular quiver of type $\tilde{\mathbb{A}}_{2n-1}$. This relation allows us to compute, independently of each other, all the cluster variables in the cluster algebra associated with a quiver $Q$ of type $\tilde{\mathbb{D}}_{n}.$ 
\vspace{0.2 cm}


\end{abstract}

\section{Introduction}

\hspace{0.7 cm} Cluster algebras were introduced by S. Fomin and A. Zelevinsky in [FZ02, FZ03]. They are a class of commutative algebras which was shown to be connected to various areas of mathematics like, for instance, combinatorics, Lie theory, Poisson geometry, Teichm\"{u}ller theory, mathematical physics and representation theory of algebras.

\hspace{0.1 cm} A cluster algebra is generated by a set of variables, called \textit{cluster variables}, obtained recursively by a combinatorial process known as \textit{mutation} starting from a set of \textit{initial} cluster variables. Explicit computation of cluster variables is difficult and has been extensively studied, see [ARS10, ADSS11, AD11, AR12, BMRRT06, BMR08].

\hspace{0.1 cm} In order to compute cluster variables, one may use \textit{friezes}, which were introduced by Coxeter [C71] and studied by Conway and Coxeter [CC73-I,II]. Various relationships are known between friezes and cluster algebras, see [CC06, M11, Pr08, ARS10, BM12, D10, AD11, ADSS12, AR12].

\hspace{0.1 cm} The present work is motivated by the use of friezes to compute cluster variables and is inspired by the results in [ARS10, AR12, Ma13 ] giving an explicit formula as a product of $2 \times 2$ matrices for all cluster variables in coefficient-free cluster algebras of types $\tilde{\mathbb{A}}$ and $\mathbb{D}$, thus explaining at the same time the Laurent phenomenon and positivity.

\hspace{0.1 cm} Our objective here is to show that similar techniques can be used for computing cluster variables in coefficient-free cluster algebras of type $\tilde{\mathbb{D}}$.

\hspace{0.1 cm} In this paper, we establish a relation between the values of a frieze of type $\tilde{\mathbb{D}}$ and some values of an $SL_{2}$-tiling $t$ associated with a particular quiver of type $\tilde{\mathbb{A}}$. For this we associate with each quiver of type $\tilde{\mathbb{D}}_{n}$ a particular quiver of type $\tilde{\mathbb{A}}_{2n-1}$. This correspondence allows us to obtain an algorithm for computing cluster variables of a cluster algebra of type $\tilde{\mathbb{D}}_{n}$ independently of each other, that is without iteration.

\hspace{0.1 cm} The article is organized as follows. In section $2$, we recall some basic notions on $SL_{2}$-tilings associated with a quiver of type $\tilde{\mathbb{A}}$. In section $3$, we set preliminaries on friezes of type $\tilde{\mathbb{D}}$ and establish a correspondence between a frieze of type $\tilde{\mathbb{D}}$ and an $SL_{2}$-tiling associated with a particular quiver of type $\tilde{\mathbb{A}}$. In section $4$, we give an algorithm to compute all the cluster variables in a cluster algebra associated with a quiver of type $\tilde{\mathbb{D}}$. 

\section{$SL_{2}$-tiling of the plane}

In this section we recall some notions on $SL_{2}$-tilings and apply them to a quiver of type $\tilde{\mathbb{A}}$ according to [ARS10].

\begin{defi}
 Let $\mathbb{K}$ be a field and $\mathbb{Z}^{2}$ the discrete plane. We call $SL_{2}$-tiling of the plane a map $t:\; \mathbb{Z}^{2}\rightarrow \mathbb{K}$ such that for any $u, v \in \mathbb{Z}$  $$ det \left( \begin{array}{cccccccccc}
t(u,v+1) &  &t(u+1,v+1)   \\
 t(u,v) &  &t(u +1,v) \\
\end{array} \right) = 1.$$\\
\end{defi}
The following example represents an $SL_{2}$-tiling with $\mathbb{K}=\mathbb{Q}$.\\

\begin{exem}

$\xymatrix@=7pt{
&&& &   &   &    &  1 &   1 & 4     & 19 &               \\
&...& &   & ...  &  &  &   1 & 2 &9   &43&       \\
&& &&   &  &  &1& 3& 14 & 67 & ...  &          \\
&&&&1&1 &1 & 1&4&19 &91&  &&&&      \\
&1&1&1&1&2&3&4&17&81&388 & ...&         \\
&1 &2&3&4 &9 & 14&19 &81& 386&1849&          \\
& & &&...&&&...&\\
}$\\
\end{exem}

The $SL_{2}$-tiling is an extension to the whole plane of the frieze introduced by Coxeter [C71] and studied by Conway and Coxeter [CC73-I, II].

\hspace{0.3 cm} Let $\Delta$ be a finite acyclic (containing no oriented cycles) quiver with $\Delta_{0}$ the set of its points and $\Delta_{1}$ the set of its arrows and $\mathbb{K}$ a field. The translation quiver $\mathbb{Z}\Delta$ associated with $\Delta$  (see [ASS-VIII.1.1]) consists of two sets: the set of points  $(\mathbb{Z}\Delta)_{0} = \mathbb{Z} \times \Delta_{0} = \left\lbrace  (k, i) | k \in \mathbb{Z}, \, i \in \Delta_{0} \right\rbrace $ and the set of arrows\\
 $(\mathbb{Z}\Delta)_{1} = \left\lbrace  (k, \alpha): (k, i)\rightarrow (k, j) | k \in \mathbb{Z}, \, \alpha: i \rightarrow j \in \Delta_{1} \right\rbrace  \cup \\ \left\lbrace  (k, \alpha'): (k, j)\rightarrow (k+1, i) | k \in \mathbb{Z}, \, \alpha: i \rightarrow j \in \Delta_{1} \right\rbrace. $ Let us define a frieze associated with the quiver $\Delta$.
 
\hspace{0.3 cm} In the translation quiver $\mathbb{Z}\Delta$ (see [ASS-VIII.1.1]), let us replace the points $(\mathbb{Z}\Delta)_{0}$ by their images under a frieze function $\mathfrak{a}$: $(\mathbb{Z}\Delta)_{0} \rightarrow \mathbb{K}$ defined for some initial values $\mathfrak{a}(0,i) \in \mathbb{K}$ as follows: $\mathfrak{a}(k,i)\mathfrak{a}(k+1,i) = 1+ \prod_{(k,i)\rightarrow (m,j)} \mathfrak{a}(m,j)$ where the product is taken over the arrows (see [ARS10-2]). The resulting translation quiver with values associated with its vertices is called a frieze.\\

In general we associate a \textit{boundary} to an $SL_{2}$-tiling.

\begin{defi}
We call boundary a bi-infinite sequence $...c_{-2}x_{-2}c_{-1}x_{-1}c_{0}x_{0}c_{1}x_{1}c_{2}x_{2}...,$ with $x_{i}\in \{ x,y \}$ and $c_{i} \in \mathbb{K}$ for all $i \in \mathbb{Z}$.
\end{defi}

\hspace{0.3 cm} Coefficients $c_{k}$ equal to one are usually omitted in the sequence representing a boundary.
 
\hspace{0.3 cm} Each boundary may be embedded into the Euclidean plane in the following way: the $x$ (or $y$) determine the horizontal (or vertical, respectively) segments of a discrete path, that is $x$ (or $y$) corresponds to a segment of the form [$(u , v)\; ,\; (u+1 , v)]$  (or $[(u , v)\; ,\; (u , v+1)]$, respectively) in the plane. The variables $c_{i}$ become thus labels of the vertices of a discrete path. A boundary is called \textit{admissible} if none of the sequences $(x_{n})_{n\leq 0}$ and $(x_{n})_{n\geq 0}$ is ultimately constant.

\hspace{0.3 cm} Given an admissible boundary $f$ embedded in the plane, let $(u,v)$ be a point in $\mathbb{Z}^{2}$. Then the \textit{word associated with} $(u,v)$ is the portion of the boundary $f$ delimited by the horizontal and vertical projections of the point $(u, v)$ on the boundary $f$.

\hspace{0.3 cm} The following example shows an embedded boundary and how to associate a word with a point in $\mathbb{Z}^{2}$.

\begin{exem}
The word associated with the point $P$ is $c_{-2}yc_{-1}xc_{0}xc_{1}yc_{2}yc_{3}xc_{4}$. \\
  
$
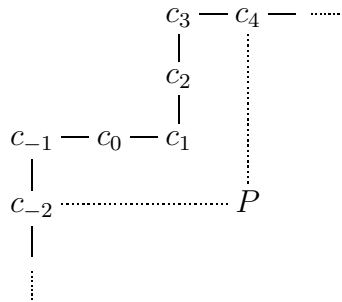
\xymatrix@=10pt{
&      &      &   &      &   c_{3}\ar@{-}[d]\ar@{-}[r]   &     c_{4}\ar@{-}[r]  &  \ar@{.}[r]    &  &         \\
&  &   &   &      & c_{2}\ar@{-}[d]     &   &       &    &  &   \\
&   &  & c_{-1}\ar@{-}[r]\ar@{-}[d] &  c_{0}\ar@{-}[r] &  c_{1}  &   &   & &  \\
&    &  & c_{-2}\ar@{-}[d] & &   &   P\ar@{.}[uuu]\ar@{.}[lll]   &      &   & \\
&    &  & \ar@{.}[d] &   &   &      &      &   & \\
&&&&&&&&&&&&&&&&\\
  }$

\vspace{1cm}

If we set $c_{i} = 1 $, for all $i$ then the word associated with the point $P$ can be written as follows: yxxyyx.\\
\end{exem}

It is always possible to construct an $SL_{2}$-tiling starting from an admissible boundary in the plane. The work [ARS10] provides a formula for a value in the tiling at the point $(u,v) \in \mathbb{Z}^{2}$. This formula is given in terms of the associated word and the following matrices:\\

$ M(a , x , b)= \left( \begin{array}{cccccc}
a &  &1   \\
0 &  &b \\
\end{array} \right) ;\quad  M(a , y , b)= \left( \begin{array}{cccccc}
b &  &0  \\
1 &  &a \\
\end{array} \right)$, \quad $a , b \in \mathbb{K}$. \\

The following theorem from [ARS10-Theorem $4$] allows us to compute the variable lying at the point $(u,v) \in \mathbb{Z}^{2}$.\\

\begin{theo}

Given an admissible boundary $f$, there exists a unique $SL_{2}$-tiling $t$ of the plane extending the embedding of the boundary into the plane. For any point $(u,v)$ below the boundary, with an associated word $b_{0}x_{1}b_{1} x_{2}...b_{n}x_{n+1}b_{n+1}$ where $n \geq 1$, $x_{i} \in \{ x,y \}$, $ b_{i} \in \mathbb{K}$, the tiling $t$ is defined by the formula\\

 $t(u,v)= \displaystyle\frac{1}{b_{1}b_{2} . . . b_{n}}(1,b_{0})\prod_{i=2}^{n}M(b_{i-1},x_{i},b_{i})\left( \begin{array}{cccccc}
1 \\
b_{n+1} \\
       
\end{array} \right)$. $\square$\\
\end{theo}

Now we give an application of the notion of $SL_{2}$-tiling to a quiver of type $\tilde{\mathbb{A}}$ according to [ARS10].\\

Consider an acyclic quiver $\Theta$ of type $\tilde{\mathbb{A}}_{n}$, $n\geq 1$, such that its points are labelled by natural numbers modulo $(n+1)$ in clockwise orientation. Let us form a word $\omega$ by associating a variable $x_{j} \in \{ x, y\}$ with arrows of $\Theta$ as follows. Let $x_{j}=x$ if the arrow is $j\rightarrow j+1$ and $x_{j}=y$ if it is $j\leftarrow j+1$. Let $\omega$ be such a word $ x_{1}x_{2}...x_{n+1}$ which encodes the orientation of arrows in the quiver $\Theta$. Let $^{\infty}\omega^{\infty}$ be the extension of $\omega$ by considering indices  $j \in \mathbb{Z}$ and putting $x_{k} = x_{j}$ if $k$ equals $j$ modulo $(n+1)$. Then $^{\infty}\omega^{\infty}$ defines an admissible boundary associated with the quiver $\Theta$.

\hspace{0.1 cm} According to [ARS10], given some initial values $\mathfrak{a}(0,i)$ the values of the tiling $t$ from theorem $2.5$ below the admissible boundary $^{\infty}\omega^{\infty}$ are values of the frieze associated with the quiver $\Theta$. These values are computed by applying the map $t$.\\

We present in the following example an $SL_{2}$-tiling for a quiver $\Theta$ of type $\tilde{\mathbb{A}}_{3}$.\\
\begin{exem}
Consider the following quiver $\Theta$ of type $\tilde{\mathbb{A}}_{3}$:\\
\resizebox{5cm}{!}{
$
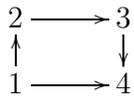
\xymatrix@=1pc{ 
 &   &2\ar[rr] &    &3\ar[d] &&&        \\
&
 \text{}   & 1\ar[rr]\ar[u]   &      & 4 & &&        \\
 &   &    &   &    &      &      &    \\
}$}\\

We have $\omega = xxxy$ and the boundary associated with $\Theta$ is $f =\, ^{\infty}(xxxy)^{\infty}$. This gives us the following $SL_{2}$-tiling:\\

$\xymatrix@=7pt{
&&& &   &   &    &   &    &   &  & &&1&              \\
&...& &   & ...  &  &  &    &  &   &1&1&1&1&       \\
&& &&   &  &  &1& 1&1 & 1 &2&3&4& &          \\
&&&&1&1 &1 & 1&2&3 &4&9  &14&19&&      \\
&1&1&1&1&2&3&4&9&14&19 &43&67& ...&         \\
&1 &2&3&4 &9 & 14&19 &43& 67&&          \\
& & &&...&&&...&\\
}$

The word corresponding to the values $9$ in the tiling is $yxxxyx$ and then, according to theorem $2.5$,\\

$9 = \left( 1,1\right)  \left( \begin{array}{cccccc}
1 &  &1 \\
0 &  &1 \\
\end{array} \right) \left( \begin{array}{cccccc}
1 &  &1 \\
0 &  &1\\
\end{array} \right) \left( \begin{array}{cccccc}
1 &  &1  \\
0 &  &1 \\
\end{array} \right) \left( \begin{array}{cccccc}
1 &  &0  \\
1 &  &1 \\
\end{array} \right)\left( \begin{array}{cccccc}
1 \\
1 \\
       
\end{array} \right).$\\
\end{exem}

\hspace{0.3 cm} An admissible boundary $f$ is said to be \textit{periodic} if it is of the form\\ $^{\infty}(c_{1}x_{1}...c_{n+1}x_{n+1}c_{1})^{\infty}$ and in this case its period is $(n+1)$. The $SL_{2}$-tiling associated with this boundary has periodicity determined by the vector $(r,s)$, where $r$ (or $s$) is the number of $x$ (or $y$, respectively) among $x_{1},...,x_{n+1}$.  We have also $r + s = n+1$.

\hspace{0.1 cm} We call the finite sequence $c_{1}x_{1}...c_{n+1}x_{n+1}c_{1}$ a generator of the boundary $f$. Note that when extending a generator $c_{1}x_{1}...c_{n+1}x_{n+1}c_{1}$ to a boundary $^{\infty}(c_{1}x_{1}...c_{n+1}x_{n+1}c_{1})^{\infty}$ we glue adjacent copies of the generator in a way that there is only one occurrence of $c_{1}$ between $x_{n+1}$ and $x_{1}$ ($...x_{n+1}c_{1}x_{1}...$). The admissible boundary $f$ associated with the quiver $\Theta$ of type $\tilde{\mathbb{A}}_{r,s}$ is $^{\infty}\omega^{\infty}$ constructed above. Its generator can be obtained by cutting the quiver $\Theta$ at one of its points and reading the resulting quiver clockwise (by reading the resulting quiver anticlockwise we obtain a generator of a boundary equivalent to $f$). This operation creates an additional point of the quiver $\Theta$, the second copy of the point at which we cut. Thus the quiver of type $\tilde{\mathbb{A}}_{r,s}$ is transformed into a quiver of type $\mathbb{A}_{r+s+1}$. 

\hspace{0.1 cm} In example $2.6$ let us cut the quiver $\tilde{\mathbb{A}}_{3,1}$ at its point $1$. We obtain in this way a quiver of type ${\mathbb{A}}_{5}$ (reading clockwise) labelled as follows: $\xymatrix @=10pt{1\ar[r]
&{2}\ar[r]
&3\ar[r]
&4
&1\ar[l]
}  $. If in this quiver of type $\mathbb{A}_{5}$ we denote by $x$ the arrows oriented toward the right and by $y$ those oriented toward the left then we get the generator $\omega = xxxy$.\\

 A quiver $\Theta$ of type $\tilde{\mathbb{A}}_{r,s}$ has $(r+s)$ possibilities of cutting, we choose any one of them. The quiver underlying $\omega$ is the quiver of type ${\mathbb{A}}_{r+s+1}$.\\

 \hspace{0.3 cm} We now recall the definition of a seed due to Fomin and Zelevinsky [FZ03-1.2].\\

\begin{defi}
 Let $\Gamma$ be a quiver (without loops or oriented two-cycles) with $(n+1)$ points and $\chi = \{u_{1},u_{2}, ..., u_{n+1}\}$ a set of variables called cluster variables, such that a variable $u_{i}$ is associated with the point $i$ (with $1\leq i \leq n+1$) of $\Gamma$. The set $\chi$ is called a cluster and the pair $(\Gamma , \chi)$ is called a seed.\\
\end{defi}

 \hspace{0.3 cm} We can obtain other seeds by mutation (see [FZ03-1.2]) starting from the seed $(\Gamma, \chi)$. The set of all cluster variables obtained by successive mutation generates an algebra over $\mathbb{Z}$ called \textit{cluster algebra} which is denoted by $\mathcal{A}(\Gamma, \chi)$. For cluster algebras, we refer to the papers [FZ02], [FZ03].
 
\begin{rem}
We say that the type of a seed $(\Gamma, \chi)$ and the cluster algebra $\mathcal{A}(\Gamma, \chi)$ coincides with the type of quiver $\Gamma$.
\end{rem}

\hspace{0.3 cm} If we associate a point $j$ of a quiver $\Theta$ of type $\tilde{\mathbb{A}}_{n}$ with the variable $u_{j}, \; j = 1, . . ., (n+1)$, these variables will correspond to the vertices of the boundary. Then the variables of the $SL_{2}$-tiling below the boundary are cluster variables of the cluster algebra $\mathcal{A}(\Theta, \{u_{1}, ..., u_{n+1}\})$. They are obtained by applying the formula of theorem $2.5$ (see [ARS10-8]; [AR12-4.4]).\\

\hspace{0.3 cm} Our aim in this paper is to compute cluster variables of a cluster algebra of type $\tilde{\mathbb{D}}_{n}$.

\section{$SL_{2}$-tiling and friezes of type $\tilde{\mathbb{D}}$ }

In this section we establish a relation between values of a frieze of type $\tilde{\mathbb{D}}$ and some values of an $SL_{2}$-tiling associated with a particular quiver of type $\tilde{\mathbb{A}}$. To this end we propose a way to associate an admissible boundary with a quiver of type $\tilde{\mathbb{D}}$.

\subsection{Preliminaries on friezes of type $\tilde{\mathbb{D}}$}

Consider a quiver $Q$ whose underlying graph is of type $\tilde{\mathbb{D}}_{n}$. We agree to label the points of $Q$ as follows:
 $$
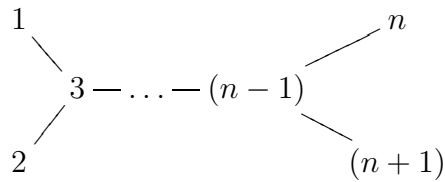
\xymatrix @=10pt
{
&1
&
&
&
&n\ar@{-}[dl]
&&\\
\text{}
&
&3 \ar@{-}[r]\ar@{-}[ul]\ar@{-}[dl]
&\ldots \ar@{-}[r]
&(n-1) 
&\text{} \\
&2
&
&
&
&(n+1)\ar@{-}[ul]
}
$$\\
 where a solid segment represents an arrow without its orientation.\\
The \textit{forks} are the full sub-quivers of $Q$ generated by the points $\{1, 2, 3\}$ and by the points $\{(n-1), n, (n+1)\}$, respectively.\\
We agree to call:\\
- \textit{fork arrows}, the arrows of each fork, \\
- \textit{joint} of each fork, the points $3$ or $n-1$, respectively\\
- and \textit{fork vertices}, the points $1$ and $2$ or $n$ and $(n+1)$, respectively.\\

\hspace{0.3 cm} If we associate with each vertex $i$ of $Q$ a variable $u_{i}$ then we get the seed $\mathcal{G}= (Q , \chi)$ whose underlying graph can be represented by the following diagram: $$
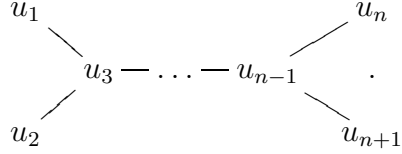
\xymatrix @=10pt
{
&u_{1}
&
&
&
&u_{n}\ar@{-}[dl]
&&\\
\text{}
&
&u_{3} \ar@{-}[r]\ar@{-}[ul]\ar@{-}[dl]
&\ldots \ar@{-}[r]
&u_{n-1} 
&\text{} .\\
&u_{2}
&
&
&
&u_{n+1}\ar@{-}[ul]
}
$$\\

\hspace{0.3 cm} Let us define a frieze on the translation quiver $\mathbb{Z}Q$ by the function\\ $\mathfrak{a}$: $(\mathbb{Z}Q)_{0} \rightarrow \mathbb{Q}(u_{1}, u_{2},...,u_{n+1})$ such that for $(k,i) \in (\mathbb{Z}Q)_{0}$ we have:\\ $\mathfrak{a}(k,i)\mathfrak{a}(k+1,i) = 1+ \prod_{(k,i)\rightarrow (m,j)} \mathfrak{a}(m,j)$ with initial variables $\mathfrak{a}(0, i) = u_{i}$ and the product taken over the arrows. Then all the values of the frieze are cluster variables of $\mathcal{A}(Q, \{u_{1}, u_{2},...,u_{n+1}\})$ (see [AD11-1.3]).\\

\hspace{0.3 cm} The study of seeds of type $\tilde{\mathbb{D}}_{n}$ can be subdivided into different cases depending on the orientation of fork arrows. For each fork we have the following three cases:\\
- The fork is composed by two arrows leaving the joint,\\
- The fork is composed by two arrows entering the joint,\\
- The fork is composed by one arrow leaving the joint and another arrow entering the joint.\\

\hspace{0.3 cm} It is well-known that for two seeds $\mathcal{G}_{1}$ and $\mathcal{G}_{2}$, which are mutation equivalent (we refer to [FZ03-8] for the notion of mutation equivalent), the cluster algebras $\mathcal{A}(\mathcal{G}_{1})$ and $\mathcal{A}(\mathcal{G}_{2})$ associated with these seeds, respectively, coincide (see [FZ03-1.2]).\\

\hspace{0.3 cm} Our aim in this paper is to compute cluster variables of a cluster algebra of type $\tilde{\mathbb{D}}_{n}$ independently of each other. The following lemma in [ASS-VII.5.2] allows us to reduce the study to one of the three different cases mentioned above.

\begin{lem}
Let $Q_{1}$ and $Q_{2}$ be two quivers having the same underlying graph $G$. If $G$ is a tree then $Q_{1}$ and $Q_{2}$ are mutation equivalent. $\square$ \\
\end{lem}

\hspace{0.3 cm} In the following, $\mathcal{G}$ denotes a seed of type $\tilde{\mathbb{D}}_{n}$ with a quiver whose each fork is composed by two arrows both entering (or leaving) the joint.\\

\hspace{0.3 cm} We denote by $F$ the frieze associated with the seed $\mathcal{G}$. We define in the following the notion of \textit{modelled quiver} which will help us to compute the cluster variables of a cluster algebra of type $\tilde{\mathbb{D}}_{n}$ without using iteration.

\begin{defi}
We call modelled quiver $\bar{{F}}$ associated with $\mathcal{G}$, the translation quiver obtained from the frieze ${F}$ as follows:\\
$1.$ by gluing in ${F}$ the arrows of each shifted copy of the forks,\\
$2.$ by multiplying the values assigned to vertices of the fork corresponding to the arrows that were glued in step $1$.  Namely the arrows obtained by gluing in step $1$ have $\mathfrak{a}(k, 1)\mathfrak{a}(k, 2)$ or $\mathfrak{a}(k, n)\mathfrak{a}(k, n+1)$ as the corresponding variables.\\

\end{defi}

We give an example of a modelled quiver associated with a seed of type $\tilde{\mathbb{D}}_{4}$.

\begin{exem}
 $\xymatrix @=10pt
{
&u_{1}\ar[dr]
&
&u_{4}\ar[dl]
&&\\
\text{For the following seed of type $\tilde{\mathbb{D}}_{4}$:}\quad
&
&u_{3} 
&
&
&\text{,}\\
&u_{2}\ar[ur]
&
&u_{5}\ar[ul]
}
$\\
the frieze ${F}$ has the form:\\
\resizebox{12cm}{!}{
$
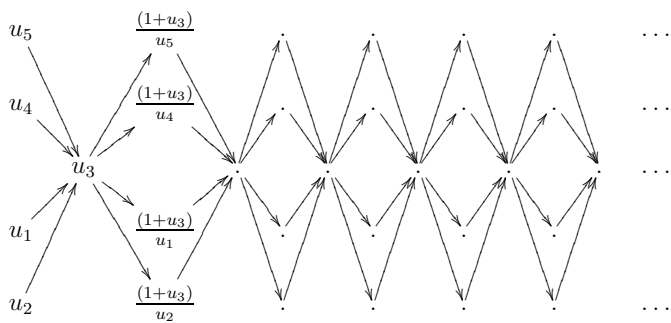
\xymatrix@=1pc{
& &   &   &    &   &    &      &      &&&&&&&&         \\
& &u_{5}\ar[ddr]&   & \frac{(1+u_{3})}{u_{5}}\ar[ddr]& & .\ar[ddr]&  & .\ar[ddr] &   & .\ar[ddr] &  &.\ar[ddr]&&\ldots& && \\
&&u_{4}\ar[dr]&   &\frac{(1+u_{3})}{u_{4}}\ar[dr] & & .\ar[dr]& &  .\ar[dr]&   &.\ar[dr]&&.\ar[dr]&      &\ldots&&\\
&&&u_{3}\ar[uur] \ar[dr]\ar[ddr]\ar[ur]&&.\ar[ddr]\ar[ur]\ar[dr]\ar[uur]& &.\ar[uur]\ar[dr]\ar[ddr]\ar[ur]& &   .\ar[uur]\ar[dr]\ar[ddr]\ar[ur]& & .\ar[dr]\ar[ddr]\ar[ur]\ar[uur] &          &.&\ldots&&&&\\
&& u_{1} \ar[ur]&&\frac{(1+u_{3})}{u_{1}}\ar[ur] &  &.\ar[ur]&  & .\ar[ur] & &  .\ar[ur] &   &.\ar[ur]&&&&      \\
&&u_{2} \ar[uur]&&\frac{(1+u_{3})}{u_{2}}\ar[uur]&&.\ar[uur] & &.\ar[uur]&&.\ar[uur]& &.\ar[uur]&&\ldots&&&&\\
 & &   &  &    &   &      &    &&&&&&&\\
}$} \\

and the modelled quiver $\bar{F}$ has the form:\\
\resizebox{12cm}{!}{
$
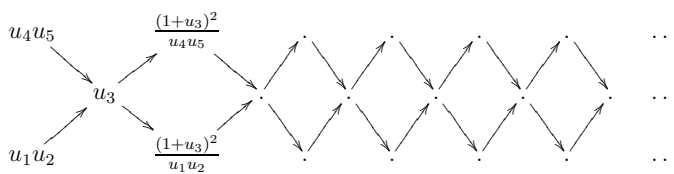
\xymatrix@=1pc{
& &   &   &    &   &    &      &      &&&&&&&&         \\
&&u_{4}u_{5}\ar[dr]&   &\frac{(1+u_{3})^{2}}{u_{4}u_{5}}\ar[dr] & & .\ar[dr]& &  .\ar[dr]&   &.\ar[dr]&&.\ar[dr]&      &\ldots&&\\
&&&u_{3} \ar[dr]\ar[ur]&&.\ar[ur]\ar[dr]& &.\ar[dr]\ar[ur]& &   .\ar[dr]\ar[ur]& & .\ar[dr]\ar[ur] &          &.&\ldots&&&&\\
&& u_{1}u_{2} \ar[ur]&&\frac{(1+u_{3})^{2}}{u_{1}u_{2}}\ar[ur] &  &.\ar[ur]&  & .\ar[ur] & &  .\ar[ur] &   &.\ar[ur]&&\ldots&&      \\
 & &   &  &    &   &      &    &&&&&&&\\
}$} \\
\end{exem}

\begin{rem}
 All the squares of the form

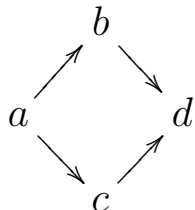
\resizebox{7cm}{!}{
$\xymatrix@=1pc{ 
 &   & &  b\ar[dr]  & &&&        \\
   &\text{}   & a\ar[ur]\ar[dr]   &      & d     & &&        \\
 &   &    &   c\ar[ur]&    &      &      &    \\
}$}\\ in $\bar{{F}}$ satisfy the relation $ad-bc = 1$ called the uni-modular rule, with $a, b, c, d \in \mathbb{Q}(u_{1}, u_{2},..., u_{n+1})$. The values on the upper and bottom extreme lines in $\bar{F}$ are products of two cluster variables (products created by the passage from $F$ to $\bar{F}$). Note that the pairs of cluster variables forming these products (in the case of both arrows of a fork entering or leaving the joint) are given by fractions whose numerators are equal and denominators coincide up to exchanging $u_{1}$ and $u_{2}$ (or $u_{n}$ and $u_{n+1}$), which appear in denominators with exponent one (see [BMR09-1]). Therefore, a product of these two variables is of the form $\displaystyle\frac{1}{(u_{1}u_{2})}\phi^{2}$ for one fork and $\displaystyle\frac{1}{(u_{n}u_{n+1})}\psi^{2}$ for the other fork, with $\phi, \psi \in \mathbb{Q}(u_{1}, u_{2},...,u_{n+1})$.\\

\hspace{0.3 cm} Then the two extreme lines in the modelled quiver $\bar{F}$ contain the sequences of variables of the form $\displaystyle\frac{1}{(u_{1}u_{2})}\phi^{2}_{k}$ for the horizontal line passing through the variable $u_{1}u_{2}$ and of the form $\displaystyle\frac{1}{(u_{n}u_{n+1})}\psi^{2}_{k}$ for the horizontal line passing through the variable $u_{n}u_{n+1}$.\\

As a consequence of the definition of a frieze, we have the following relations for all $k \in \mathbb{N}$ : $$(1)\hspace{4cm} [\mathfrak{a}(k,3) + 1]^{2} = \displaystyle\frac{1}{(u_{1}u_{2})^{2}}\phi^{2}_{k}\phi^{2}_{k+1}, \hspace{3 cm} $$
 $$(2)\hspace{2.7cm} [\mathfrak{a}(k,n-1) + 1]^{2} = \displaystyle\frac{1}{(u_{n}u_{n+1})^{2}}\psi^{2}_{k}\psi^{2}_{k-1}. \hspace{3 cm} $$\\

\end{rem}

\subsection{Correspondence between a frieze of type $\tilde{\mathbb{D}}$ and an $SL_{2}$-tiling associated with a particular quiver of type $\tilde{\mathbb{A}}$}

Let $Q$ be a quiver of type $\tilde{\mathbb{D}}_{n}$ and \textbf{$\Sigma$} be the full sub-quiver of $Q$ generated by all points except the points $2$, $(n+1)$. We agree to draw $\Sigma$ from the left to the right in such a way that the vertices appear in increasing order. Note that $\Sigma$ is a quiver of type $\mathbb{A}_{n-1}$. This leads to a way to associate a quiver of type $\tilde{\mathbb{A}}_{2n-1}$ with a quiver of type $\tilde{\mathbb{D}}_{n}$, which, in turn, will allow us to associate an admissible boundary with a quiver of type $\tilde{\mathbb{D}}_{n}$.

\hspace{0.3 cm} For a quiver $\Lambda$ of type $\mathbb{A}_{n}$ (drawn from left to right), let us denote by $^{t}\Lambda$ the transpose of $\Lambda$, that is, the quiver obtained by redrawing $\Lambda$ from right to left.\\

\hspace{0.3 cm}  With a quiver $Q$ of type $\tilde{\mathbb{D}}_{n}$ we associate the quiver $Q'$ of the form\\
\resizebox{6cm}{!}
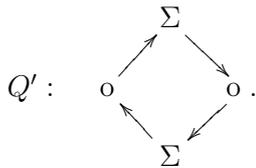
{
$$\xymatrix@=1pc{ 
 &   & &  \Sigma\ar[dr]  & &&&        \\
   &Q':   & \text{o}\ar[ur]   &      & \text{o}\; . \ar[dl]   & &&        \\
 &   &    &   \Sigma\ar[ul]&    &      &      &    \\
}$$}\\

\hspace{0.3 cm} We know that the boundary $f =\, ^{\infty}\omega^{\infty}$ associated with a quiver $\Theta$ of type $\tilde{\mathbb{A}}_{r,s}$ is periodic and its period $(r+s)$ corresponds to the length of a generator $\omega$. 

\hspace{0.3 cm} For the quiver $Q'$ of type $\tilde{\mathbb{A}}_{2n-1}$ there is a choice of $2n$ different generators for its associated boundary $f'$: one can cut at any of $2n$ points of $Q'$ to obtain a generator. We are going to use only one of these generators to construct the admissible boundary $\tilde{f}$ associated with a quiver of type $\tilde{\mathbb{D}}_{n}$.

\hspace{0.3 cm} The generator $\omega_{1}$ of the boundary $f'$ is obtained by cutting the quiver $Q'$ at the point $1$ of its arrow $\xymatrix @=10pt{\text{o}\ar[r]
&{1}}  $
 and reading the quiver clockwise. Consider also a generator $\omega_{2}$ of an boundary equivalent to $f'$: $\omega_{2}$ is obtained by cutting $Q'$ at the point $n$ of its arrow $\xymatrix @=10pt{\text{o}\ar[r]
&{n}}$ and reading the quiver anticlockwise (see example $3.6$). We get in this way two different generators associated with the quiver $Q'$, that is, $\omega_{1}$ with the underlying quiver $\xymatrix @=10pt{{{\Sigma}}\ar[r]
&\text{o}\ar[r]
&^{t}{\Sigma}\ar[r]
&\text{o}\ar[r]
&1}  $
and $\omega_{2}$ with the underlying quiver $\xymatrix @=10pt{{n}
&\text{o}\ar[l]
&^{t}{\Sigma}\ar[l]
&\text{o}\ar[l]
&{\Sigma}\ar[l]}.$ 

\hspace{0.3 cm} Let us denote by $\bar{\omega}$ the part of $Q$ whose underlying quiver is $\Sigma$. Then we have $\omega_{1} = \bar{\omega}xx ^{t}\bar{\omega}xx$ and $\omega_{2} = yy ^{t}\bar{\omega}yy\bar{\omega}$.

\hspace{0.3 cm} Now we are in a position to associate an admissible boundary with a quiver of type $\tilde{\mathbb{D}}_{n}$ whose each fork consists of two arrows entering or leaving the joint.\\

\begin{defi}
Let $Q$ be a quiver of type $\tilde{\mathbb{D}}_{n}$ whose each fork consists of two arrows both entering or both leaving the joint. The admissible boundary $\tilde{f}$ associated with $Q$ is of the form $\tilde{f} =\,    ^{\infty}(\omega_{1})\bar{\omega}(\omega_{2})^{\infty}$, where $\omega_{1}$ is repeated infinitely many times on the left and $\omega_{2}$ is repeated infinitely many times on the right.
\end{defi}

\hspace{0.3 cm} The boundary $\tilde{f}$ is obtained by gluing to $\bar{\omega}$ the generators ${\omega}_{1}$ periodically on the left and $\omega_{2}$ periodically on the right. Note that $\bar{\omega}$ is also contained in $\omega_{1}$ and $\omega_{2}$. However the decomposition of the boundary of the form $\tilde{f} =\,    ^{\infty}(\omega_{1})\bar{\omega}(\omega_{2})^{\infty}$ defines a distinguished occurrence of $\bar{\omega}$ in $\tilde{f}$. It is this occurrence of $\bar{\omega}$ which we call the \textit{root} of $Q$ in $\tilde{f}$. \\

\hspace{0.3 cm} We give an example showing how to associate an admissible boundary with a quiver $Q$ of type $\tilde{\mathbb{D}}_{n}$.
\begin{exem}
 $\xymatrix @=10pt
{
&{1}\ar[dr]
&
&{4}\ar[dl]
&&\\
\text{Consider the following quiver $Q$ of type $\tilde{\mathbb{D}}_{4}$:}\quad
&
&{3} 
&
&
&\text{,}\\
&{2}\ar[ur]
&
&{5}\ar[ul]
}
$\\
$
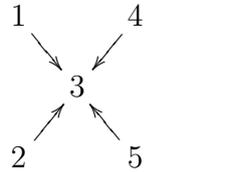
\xymatrix @=10pt{
&
&
&\text{}\\
& \text{we have} \; \; \textbf{$\Sigma$} :\;{1}\ar[r]
&{3}
&4\ar[l]
&&\text{}\\  
&
&
&
&
}  $
$\xymatrix @=10pt{
&
&
&\text{}\\
& \text{and} \; \; \textbf{$^{t}\Sigma$} :\;{4}\ar[r]
&{3}
&1.\ar[l]
&&\text{}\\  
&
&
&
&
}  $\\
 The quiver $Q'$ of type $\tilde{\mathbb{A}}_{7}$ associated with $Q$ is:\\
\resizebox{6cm}{!}{
$
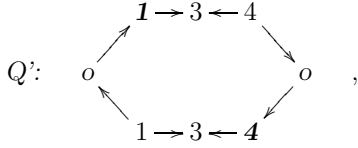
\xymatrix@=1pc{ 
 &   & & \textbf{1}\ar[r]  &3 &4\ar[dr]\ar[l]&&        \\
   &\text{Q':}   & \text{o}\ar[ur]   &      &    & &\text{o}\ar[dl]&,&        \\
 &   &    &   1\ar[ul]\ar[r]&3    &  \textbf{4}\ar[l]    &      &    \\
}$}\\
This gives us two generators $\omega_{1}$ and $\omega_{2}$ associated with $Q'$ such that their underlying quivers are respectively\\ $\xymatrix @=10pt{
&\text{${\omega}_{1}$} :
&1\ar[r]
&3
&4\ar[l]\ar[r]
&\text{o}\ar[r]
&4\ar[r]
&3
&1\ar[r]\ar[l]
&\text{o}\ar[r]
&1
&
}  $
and \\$\xymatrix @=10pt{
&\text{${\omega}_{2}$} :
&4
&\text{o}\ar[l]
&4\ar[l]\ar[r]
&{3}
&1\ar[l]
&\text{o}\ar[l]
&1\ar[r]\ar[l]
&3
&4.\ar[l]
&
}$\\

\hspace{0.3 cm} The generators $\omega_{1}$ and $\omega_{2}$ are obtained by cutting the quiver $Q'$ of type $\tilde{\mathbb{A}}_{7}$ at its points $1$ and $4$ represented in bold characters on the picture of $Q'$, respectively.\\

\hspace{0.3 cm} The boundary $\tilde{f}$ associated with $Q$ is: $$\tilde{f} =\,    ^{\infty}(\omega_{1})\bar{\omega}(\omega_{2})^{\infty} = \; ^{\infty}(xyxxxyxx)(xy)(yyxyyyxy)^{\infty}.$$
\end{exem}

\hspace{0.3 cm} From an admissible boundary and given initial values it is always possible to construct an $SL_{2}$-tiling $t$, therefore we construct an $SL_{2}$-tiling $t$ from the boundary $\tilde{f}$ (with all initial values equal one).\\

\hspace{0.3 cm} We now recall the notion of \textit{ray} in a tiling $t$ from [ARS10-6.1]. This notion allows us to find the values of a frieze of type $\tilde{\mathbb{D}}_{n}$ among values in the tiling $t$ below the boundary $\tilde{f}$.\\

\hspace{0.3 cm} Given a mapping  $t:\; \mathbb{Z}^{2}\rightarrow \mathbb{K}$, a point $M \in \,\mathbb{Z}^{2}$ and a nonzero vector $V \in \,
\mathbb{Z}^{2}$, we consider the sequence $a_{n} = t(M+nV), n \in \mathbb{N}$. Such a sequence will be called a \textit{ray} associated with $t$. We call $M$ the \textit{origin} of the ray and $V$ its \textit{directing vector}. The ray is \textit{horizontal} if $V = (1,0)$, \textit{vertical} if $V = (0,1)$ and \textit{diagonal} if $V = (1,1)$.\\

\begin{exem}
Let $\mathbb{K} = \mathbb{Q}$ and consider the quiver $Q$ of type $\tilde{\mathbb{D}}_{4}$ of example $3.6$.\\
The boundary $\tilde{f}$ associated with $Q$ is $\tilde{f} =\; ^{\infty}(xyxxxyxx)(xy)(yyxyyyxy)^{\infty}$. This gives us the following $SL_{2}$-tiling below the boundary $\tilde{f}$: \\

$\xymatrix@=7pt{
&&&&&&&&&&&&1&3&\\
&&&&&&&&&&&1&1&4&&\\
&&&&&&&&&&&1&2&9&&&&\\
&&&&&&&&&&&1&3&14&&&&\\
&&&&&& &   &   &    &  1 &   1 & 4     & 19 &               \\
&&&&& &   &  &  &  &   1 & 2 &9   &43&       \\
&&&&& &&   &  &  &1& \textbf{3}& \textbf{14} & \textbf{67} & ...  &          \\
&&&&&&&1&1 &1 & 1&4&19 &91&  &&&&      \\
&&&&1&1&1&1&\textbf{2}&3&4&\textbf{17}&81&388 & ...&         \\
&1&1&1&1 &2&3&4 & \textbf{9}& 14&19 &81&\textbf{386} &1849&          \\
&1&2&3&4&9&14&19&\textbf{43}&67&91&3881&1849&\textbf{8857}&&&&&\\
&&&&& & &&...&&&...&\\
}$\\

\hspace{0.3 cm} The values represented in bold characters in the tiling are the vertical (diagonal , horizontal) rays with origin $2$ ( $17$, $3$ respectively).\\
\end{exem}

\hspace{0.3 cm} Now we associate variables with the points of $\Sigma$ as follows:\\
$1.$ with a vertex labelled by $i \in \{1, n \}$ we associate the products $u_{1}u_{2}$, $u_{n}u_{n+1}$ respectively,\\
$2.$ with a vertex labelled by $i \notin \{1, n \}$ we associate the variable $u_{i}$.

\begin{rem}
Starting with the quiver $\Sigma$ with associated variables as above, we can reconstruct the modelled quiver $\bar{F}$ associated with the seed $\mathcal{G}$ by using the uni-modular rule and formulas $(1)$ and $(2)$ of remark $3.4$.  
\end{rem}

\hspace{0.3 cm} If we associate also the vertex of $Q'$ labelled $i = \text{o}$ with the variable $u_{0}$ then the result of this association is a seed which will be denoted by $\mathcal{G}'$. Its underlying graph with variables is:\\
\resizebox{9cm}{!}{
$
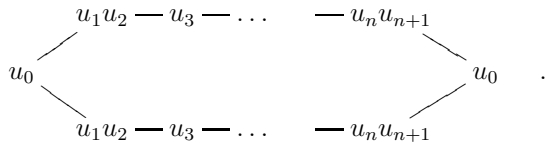
\xymatrix@=1pc{ 
 &   & & u_{1}u_{2}\ar@{-}[r]  &u_{3}\ar@{-}[r]&\ldots & \ar@{-} &u_{n}u_{n+1}\ar@{-}[dr]\ar@{-}[l]&&        \\
   &   & u_{0}\ar@{-}[ur]   & &  &   &    & &u_{0}\ar@{-}[dl]&.&        \\
 &   &    &   u_{1}u_{2}\ar@{-}[ul]\ar@{-}[r]&u_{3}\ar@{-}[r]&\ldots &\ar@{-}&  u_{n}u_{n+1}\ar@{-}[l]    &      &    \\
}$}\\

\hspace{0.3 cm} Since $u_{0}$ is not a cluster variable in $\mathcal{G}$, in the rest of this paper we evaluate $u_{0} = 1$.

\hspace{0.3 cm} Note that the generators $\bar{\omega}, \omega_{1}$ and $\omega_{2}$ are obtained from $Q'$ and thus our association of variables with the points of $Q'$ naturally associates variables with the vertices of $\bar{\omega}, \omega_{1}$ and $\omega_{2}$. 

\hspace{0.3 cm} From now on, we denote by $\bar{\omega}, \omega_{1}$ and $\omega_{2}$ the corresponding generators with variables.

\hspace{0.3 cm} We denote by $\tilde{f}_{0}$ the admissible boundary $\tilde{f}$ with variables and the following theorem establishes a link between the frieze $F$ and the tiling associated with $\tilde{f}_{0}$.

\begin{theo}
Let $n\geq 4$ be an integer, $\mathcal{G}$ a seed of type $\tilde{\mathbb{D}}_{n}$ and $\bar{F}$ the modelled quiver associated with $\mathcal{G}$. Then the horizontal lines in $\bar{F}$ are the diagonal rays in the $SL_{2}$-tiling associated with $\tilde{f}_{0}$ such that the origin of each ray is a vertex of the root $\bar{\omega}$ in $\tilde{f}_{0} =\,    ^{\infty}(\omega_{1})\bar{\omega}(\omega_{2})^{\infty}$.
\end{theo}

\textbf{Proof}

\hspace{0.3 cm} Consider the admissible boundary $$\tilde{f}_{0} =\,    ^{\infty}(\omega_{1})\bar{\omega}(\omega_{2})^{\infty} =\; ^{\infty}(\bar{\omega}xx\,^{t}\bar{\omega}xx)\bar{\omega}(yy\, ^{t}\bar{\omega}yy\bar{\omega})^{\infty}. $$

\hspace{0.3 cm} By doing a translation of parentheses toward extremities, we can rewrite $\tilde{f}_{0}$ as follows: $$\tilde{f}_{0} =\; ^{\infty}(xx\bar{\omega}xx\,^{t}\bar{\omega})xx\bar{\omega}yy( ^{t}\bar{\omega}yy\bar{\omega}yy)^{\infty}.$$ 

\hspace{0.3 cm} Let us denote $s =\, ^{\infty}(xx\bar{\omega}xx\,^{t}\bar{\omega})$ and $s' =\, ( ^{t}\bar{\omega}yy\bar{\omega}yy)^{\infty}$. The following figure illustrates the admissible boundary $\tilde{f}_{0}$ in its form $\tilde{f}_{0} = sxx\bar{\omega}yys'$.

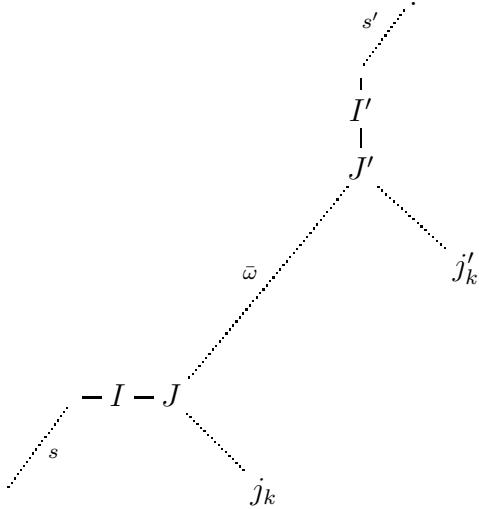
$$\xymatrix@=7pt{
&&& &   &   &   &  &&.&      &       &&&&        \\
 &   &   &  &  &    &     &   &   &  &&&&     \\
 &&  &  &  &&   & &\ar@{.}[ruu]^{s'} &  &   &          \\
&&  & & && & & I' \ar@{-}[u]& &&&&&&& &         \\
&&& & & &   & & J'\ar@{.}[rrdd]\ar@{-}[u]  &     &&&&&&     \\
& &  & &  &  &   &   &   &&&&&&       \\
& &  & &  &  &   &   &   & &j'_{k}&&      \\
& &  & & &  &   &   &          \\
&&&&&&&&&&&&&&&&&&\\
&&\ar@{.}[dlld]^{s} &I\ar@{-}[l]&J\ar@{.}[uuuuurrrr]^{\bar{\omega}}\ar@{.}[rrdd]\ar@{-}[l]&&&&&&\\
 && &   &  &    &   & &&  &&&&&&&  &     \\
&&&&&&j_{k}&&&&&&&&&&&&&}$$\\

\hspace{0.3 cm} Our aim is to prove that the horizontal lines in the modelled quiver $\bar{F}$ coincide with the diagonal rays with origins at the vertices of the root $\bar{\omega}$ in $\tilde{f}_{0}$, that is, with the oblique band delimited by the diagonal rays with origins $J$ and $J'$ (inclusively) in the figure above. To this end, we first calculate variables lying at the points $j_{k}$ and $j'_{k}$ and prove that they have the form of the variables on the top and bottom horizontal lines of $\bar{F}$, see remark $3.4$.\\

\hspace{0.3 cm} Consider the transpose of $s$, which is $^{t}s =\, (\bar{\omega}yy\, ^{t}\bar{\omega}yy)^{\infty}$. By a translation of parentheses toward the right, we can rewrite $^{t}s$ as follows: $$^{t}s =\, \bar{\omega}yy( ^{t}\bar{\omega}yy\bar{\omega}yy)^{\infty} = \bar{\omega}yys'. $$

\hspace{0.3 cm} This allows us to rewrite the boundary in the form:  $\tilde{f}_{0} = sxx\, ^{t}s$. Then the embedding of $\tilde{f}_{0}$ in the plane gives the following scheme:\\
$
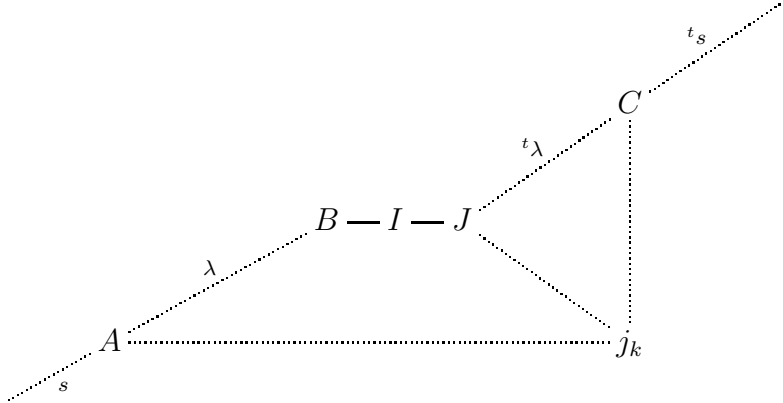
\xymatrix@=1pc{
&&&&&&&&&&&\\
&&&&&&&&&&&&&\\
&&&&&&&&&&&&&&\\
&&  & &&& && &&   &           \\
&&& &&& & & &   & & C\ar@{.}[dddd]\ar@{.}[ururr]^{^{t}s} &          \\
&&&& &&&  & &  & &   &   &          \\
&&&& &  & B\ar@{-}[r]  &I\ar@{-}[r] & J \ar@{.}[urrur]^{^{t}\lambda} &   &          \\
&&&&&& &  & &  &  &   &   &          \\
&&A\ar@{.}[rrrrrrrrr]\ar@{.}[rrrruu]^{\lambda}\ar@{.}[dll]^{s}& &&& &&&&&j_{k}\ar@{.}[llluu]&&\\
 &&&& &   &  &    &   & &&    &     \\
 &&&&&&&&&&&&&\\
 &&&&&&&&&&&&&\\
}$\\

\hspace{0.3 cm} Let us write $\lambda$ for the walk on the boundary from $A$ to $B$ and $^{t}\lambda$ for that from $J$ to $C$. Then the word associated with $j_{k}$ is of the form $\lambda xx\, ^{t}\lambda$.

\hspace{0.3 cm} By construction of $\tilde{f}_{0}$, the end point of the walk $\lambda$ is the variable $u_{1}u_{2}$, and thus $u_{1}u_{2}$ is also the starting vertex of $^{t}\lambda$. Computing the tiling function $t$ from theorem $2.5$ for the point $j_{k}$, gives:

$$t(j_{k}) = \displaystyle\frac{1}{(u_{1}u_{2})^{2}\gamma^{2}}\left( a,b\right)  \left( \begin{array}{cccccc}
u_{1}u_{2} &  &1\\
0 &  &1 \\
\end{array} \right) \left( \begin{array}{cccccc}
1 &  &1  \\
0 &  &u_{1}u_{2}\\
\end{array} \right)\left( \begin{array}{cccccc}
a \\
b \\
       
\end{array} \right),$$\\
where $\gamma$ is a product of the variables in $\lambda$ except the variable $u_{1}u_{2}$, which are also by transposition in $^{t}\lambda$, and $a, b \in \mathbb{Q}(u_{1}, u_{2},...,u_{n+1})$ arise from the product of matrices corresponding to the walk $\lambda$.

\hspace{0.3 cm} The computation gives 
$$(3)\hspace{2 cm} t(j_{k}) = \displaystyle\frac{1}{u_{1}u_{2}}\left[\frac{a+b}{\gamma}\right]^{2}.\hspace{3 cm}$$ 

\hspace{0.3 cm} Thus the value lying at the point $j_{k}$ on the diagonal ray with origin $J$ is a perfect square divided by $u_{1}u_{2}$.

\hspace{0.3 cm} To show that the value lying at a point $j'_{k}$ on the diagonal ray with origin $J'$ is a perfect square divided by $u_{n}u_{n+1},$ we perform a similar calculation using the following form of the boundary $\tilde{f}_{0} =\, ^{t}s'yys'$. This form is obtained by noticing that the transpose of $s'$ is $^{t}s' =\, ^{\infty}(xx\,^{t}\bar{\omega}xx\bar{\omega}) $ and by a translation of parentheses toward the left, we can rewrite $^{t}s'$ as follows:\\ $$^{t}s' =\, ^{\infty}(xx\bar{\omega}xx\,^{t}\bar{\omega})xx\bar{\omega} = sxx\bar{\omega}.$$\\ 

\hspace{0.3 cm} Since the tiling $t$ (below $\tilde{f}_{0}$) satisfies the uni-modular rule and since $\bar{\omega}$ with variables coincides with the quiver $\Sigma$ with variables then by virtue of remark $3.8$ it only remains to prove that the variables on the two diagonal rays with origins $J$ and $J'$ satisfy relations of the form $(1)$ and $(2)$ from remark $3.4$.\\

\hspace{0.3 cm} Consider the following scheme which represents four adjacent points $ j_{k}, j_{k+1}, d_{k}, i_{k+1}$ forming a square in the tiling below $\tilde{f}_{0}$.

$
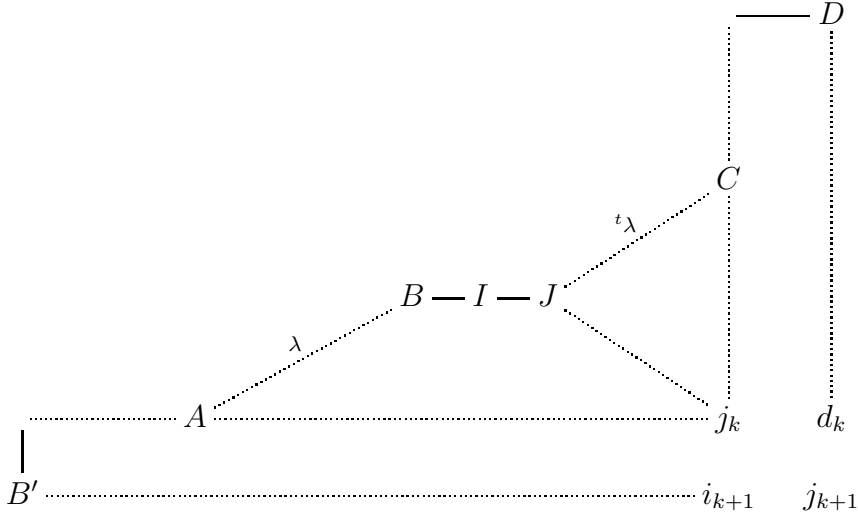
\xymatrix@=1pc{
&&&&&&&&&&&&&&&&\\
&&&&&&&&&&&&&&\ar@{-}[r]&D&&&&&&\\
&&&&&&&&&&&&&&&&&&&&&&&&\\
&&&&&&  & &&& && &&   & &&&&&&&&&    \\
&&& &&&&&&& & &    & & C\ar@{.}[dddd]\ar@{.}[uuu] &          \\
& &&&  &&&&& &  & &   &   &          \\
&&& &&&&&  & B\ar@{-}[r]  &I\ar@{-}[r] & J \ar@{.}[urrur]^{^{t}\lambda} &   &          \\
&&& &  &&&&& &  &  &   &   &          \\
&&\ar@{-}[d]&&&A\ar@{.}[lll]\ar@{.}[rrrrrrrrr]\ar@{.}[rrrruu]^{\lambda}& &&& &&&&&j_{k}\ar@{.}[llluu]&d_{k}\ar@{.}[uuuuuuu]&\\
 &&B'&&&& &   &  & &   & && &i_{k+1}\ar@{.}[llllllllllll] &j_{k+1}&&&&&&&&&&    \\
&&&&&&&&&&&&&&&&&&&&&\\
}$\\

\hspace{0.3 cm} The words associated with the points $ j_{k}, j_{k+1}, d_{k}$ and $i_{k+1}$ are respectively of the form $\lambda xx ^{t}\lambda, \mu xx ^{t}\mu, \lambda xx ^{t}\mu $ and $\mu xx ^{t}\lambda$ where $\lambda$ is the walk on the boundary from $A$ to $B$ and $\mu$ is that from $B'$ to $B$.\\

\hspace{0.3 cm} By definition of $\tilde{f}_{0}$, the point $d_{k}$ is on the diagonal ray with origin $u_{3}$. We need to prove that each value $t(d_{k})$ on this diagonal satisfies the relation: $$ (4)\hspace{2 cm} [t(d_{k}) + 1]^{2} = t(j_{k})t(j_{k+1}). \hspace{3 cm}$$

\hspace{0.3 cm} To this end, we compute the values lying at the points $ j_{k}, j_{k+1}, d_{k}$ and $i_{k+1}$. \\

\hspace{0.3 cm} Due to $(3)$, we have:
\begin{eqnarray}
t(j_{k})&=& \displaystyle\frac{1}{u_{1}u_{2}}\left[\frac{a+b}{\gamma}\right]^{2}, \nonumber \\
t(j_{k+1}) &=& \displaystyle\frac{1}{u_{1}u_{2}}\left[\frac{\tilde{a}+\tilde{b}}{\tilde{\gamma}}\right]^{2}. \nonumber 
\end{eqnarray}

For the two remaining variables, the application of the map $t$ gives:

\begin{eqnarray}
t(d_{k}) &=& \displaystyle\frac{1}{(u_{1}u_{2})^{2}\gamma \tilde{\gamma}}\left( a,b\right)  \left( \begin{array}{cccccc}
u_{1}u_{2} &  &1\\
0 &  &1 \\
\end{array} \right) \left( \begin{array}{cccccc}
1 &  &1  \\
0 &  &u_{1}u_{2}\\
\end{array} \right)\left( \begin{array}{cccccc}
\tilde{a} \\
\tilde{b} \\  
\end{array} \right) \nonumber \\
&=& \displaystyle\frac{1}{u_{1}u_{2}\gamma\tilde{\gamma} }\left[a(\tilde{a} + 2\tilde{b}) + b\tilde{b}\right], \nonumber \\
t(i_{k+1}) &=& \displaystyle\frac{1}{(u_{1}u_{2})^{2}\gamma \tilde{\gamma}}\left( \tilde{a},\tilde{b}\right)  \left( \begin{array}{cccccc}
u_{1}u_{2} &  &1\\
0 &  &1 \\
\end{array} \right) \left( \begin{array}{cccccc}
1 &  &1  \\
0 &  &u_{1}u_{2}\\
\end{array} \right)\left( \begin{array}{cccccc}
a \\
b \\  
\end{array} \right) \nonumber \\
&=& \displaystyle\frac{1}{u_{1}u_{2}\gamma\tilde{\gamma} }\left[\tilde{a}(a + 2b) + b\tilde{b}\right]. \nonumber
\end{eqnarray} 

\hspace{0.3 cm} Then we have  $t(j_{k})t(j_{k+1}) - t(d_{k})t(i_{k+1}) = \displaystyle\frac{(\tilde{a}b-a\tilde{b})^{2}}{(u_{1}u_{2}\gamma \tilde{\gamma})^{2}}$. According to the uni-modular rule in the tiling  $t(j_{k})t(j_{k+1}) - t(d_{k})t(i_{k+1}) = 1$, therefore we get $\tilde{a}b-a\tilde{b} = u_{1}u_{2}\gamma\tilde{\gamma} .$\\

\hspace{0.3 cm} Let us compute now $ t(d_{k}) + 1.$\\

\begin{eqnarray}
 t(d_{k}) + 1 &=& \displaystyle\frac{1}{u_{1}u_{2}\gamma \tilde{\gamma}}\left[a (\tilde{a} + 2\tilde{b}) + b\tilde{b}\right] + 1 \nonumber \\
&=& \displaystyle\frac{1}{u_{1}u_{2}\gamma \tilde{\gamma} }\left[a(\tilde{a} + 2\tilde{b}) + b\tilde{b} + u_{1}u_{2}\gamma \tilde{\gamma}\right] \nonumber \\
&=& \displaystyle\frac{1}{u_{1}u_{2}\gamma \tilde{\gamma}}\left[a(\tilde{a} + 2\tilde{b}) + b\tilde{b} + (\tilde{a}b - a\tilde{b})\right] \nonumber \\
&=& \displaystyle\frac{(a + b) ( \tilde{a}+\tilde{b} )}{u_{1}u_{2}\gamma \tilde{\gamma} }. \nonumber  
\end{eqnarray}

\hspace{0.3 cm} Then we have
\begin{eqnarray}
 [t(d_{k}) + 1 ]^{2} &=& \displaystyle\frac{(a + b)^{2}}{u_{1}u_{2}\gamma^{2}} \times \displaystyle\frac{(\tilde{a} + \tilde{b})^{2}}{u_{1}u_{2} \tilde{\gamma}^{2}} \nonumber \\
 &=& t(j_{k})t(j_{k+1}). \nonumber
 \end{eqnarray}

\hspace{0.3 cm} We prove in the same way that the sequence of the points $d'_{k}$ lying on the diagonal ray with origin $u_{n-1}$ satisfies the relation:   $$(5) \hspace{3 cm} [(d'_{k}) + 1 ]^{2} =  t(j'_{k})t(j'_{k-1}). \hspace{3 cm} $$

\hspace{0.3 cm} All the squares in the oblique band delimited by the diagonal rays with origins $J$ and $J'$ satisfy the uni-modular rule due to properties of the $SL_{2}$-tiling $t$. This and relations $(4)$ and $(5)$, by virtue of remark $3.8$, allows us to conclude that the diagonal rays with origins at the vertices of the root $\bar{\omega}$ of the quiver $Q$ in $\tilde{f}_{0}$ coincide with horizontal lines of the modelled quiver $\bar{F}$ associated with the seed $\mathcal{G} = (Q, \chi).$ $\square$\\

\begin{rem}
We proved theorem $3.9$ for a seed $\mathcal{G}$ whose associated quiver has forks composed by two arrows entering (or leaving) the joint. For a seed $\mathcal{G}_{2}$ whose associated quiver $Q_{2}$ has at least one fork consisting of one arrow entering and one arrow leaving the joint there are two possibilities.

\hspace{0.3 cm} One possibility is to perform a mutation on one of the fork vertices and thus reduce this case to the one considered in the proof. 

\hspace{0.3 cm} Another possibility is to work directly with the seed $\mathcal{G}_{2} = (Q_{2}, \chi)$. In this case the given proof can easily be modified, namely one has to associate with the vertex labelled by $i=1$ (or $i = n$) in $\Sigma_{2}$ the variable $\displaystyle\left( \frac{u_{2}(1+u_{3})}{u_{1}}\right)$ (or $\displaystyle\left( \frac{u_{n+1}(1+u_{n-1})}{u_{n}}\right)$) and to consider the following underlying graph for $\Sigma_{2}$ with variables:\\
$\xymatrix @=10pt
{&\displaystyle\left( \frac{u_{2}(1+u_{3})}{u_{1}}\right)  \ar@{-}[r]
&u_{3} \ar@{-}[r]
&u_{4} \ldots
&u_{n-1} \ar@{-}[l]
&u_{n}u_{n+1} \ar@{-}[l]
& 
}
$ \\
or  $\xymatrix @=10pt
{&u_{1}u_{2} \ar@{-}[r]
&u_{3} \ar@{-}[r]
&u_{4} \ldots
&u_{n-1} \ar@{-}[l]
&\displaystyle\left( \frac{u_{n+1}(1+u_{n-1})}{u_{n}}\right)  \ar@{-}[l]}
$, respectively .
\end{rem}

\section{Computation of cluster variables: case $\tilde{\mathbb{D}}_{n}$} 
Let $\mathcal{G} = (Q, \{u_{1}, ..., u_{n+1}\})$ be a seed with $Q$ of type $\tilde{\mathbb{D}}_{n}$ with forks composed by two arrows both entering or both leaving the joint.

\hspace{0.3 cm} In this section we compute cluster variables of the cluster algebra $\mathcal{A}(\mathcal{G})$ of type $\tilde{\mathbb{D}}_{n}$ by an explicit formula using matrix product (from theorem $2.5$) and by using signed continuant polynomials.\\

\hspace{0.3 cm} There are two kinds of cluster variables of a cluster algebra $\mathcal{A}(\mathcal{G})$ of type $ \tilde{\mathbb{D}}_{n}$. Cluster variables of $\mathcal{A}(\mathcal{G})$ either are transjective and correspond to those lying in the translation quiver $\mathbb{Z}Q$ as embedded in the Auslander-Reiten quiver of the cluster category or are nontransjective and correspond to those lying on the tubes (see [SS-X.1]).\\

\subsection{Computation of transjective cluster variables of a cluster algebra of type $\tilde{\mathbb{D}}_{n}$}

\hspace{0.3 cm} It is well-known, see [ARS10] and also [AD11] that for $\mathfrak{a}(0, i) = u_{i}, \, i \in (\tilde{\mathbb{D}}_{n})_{0}$, the values contained in the frieze $F$ give all transjective cluster variables of $\mathcal{A}(\mathcal{G})$ . Therefore the values of the frieze at the points of the modelled quiver $\bar{F}$ associated with the seed $\mathcal{G}$ are either transjective cluster variables of $\mathcal{A}(\mathcal{G})$ or products of two transjective cluster variables of $\mathcal{A}(\mathcal{G})$. 

\hspace{0.3 cm} According to theorem $3.9$ the modelled quiver $\bar{F}$ associated with the seed $\mathcal{G}$ is contained in the $SL_{2}$-tiling below the boundary $\tilde{f}_{0}$.

\hspace{0.3 cm} Therefore to compute transjective cluster variables of the cluster algebra $\mathcal{A}(\mathcal{G})$ of type $\tilde{\mathbb{D}}_{n}$, we proceed as follows:

\hspace{0.3 cm} In the $SL_{2}$-tiling $t$, the variables lying on diagonal rays with origins $u_{i}, \\ i = 3, 4,...,n-1$, respectively, which are vertices of the root $\bar{\omega}$ in $\tilde{f}_{0}$, are the transjective cluster variables. Their computation is obtained by applying the formula of theorem $2.5$.

\hspace{0.3 cm} The variables lying on the diagonal rays with origins $u_{1}u_{2}$ or $u_{n}u_{n+1}$, respectively, which are extreme vertices of the root $\bar{\omega}$ in $\tilde{f}_{0}$, are products of two transjective cluster variables (products created by the passage from $F$ to $\bar{F}$). Recall that the pairs of cluster variables forming these products are given by fractions whose numerators are equal and denominators coincide up to exchanging $u_{1}$ and $u_{2}$ (or $u_{n}$ and $u_{n+1}$), which appear in denominators with exponent one (see [BMR09-1]). Therefore, a product of these two transjective variables is a perfect square divided by $u_{1}u_{2}$ (or $u_{n}u_{n+1}$, respectively).

\begin{Cor}
Each value $t(u,v)$ lying on the diagonal rays with origins $u_{1}u_{2}$ or $u_{n}u_{n+1}$ gives rise to two transjective cluster variables $U$ and $V$ as follows: \\

$t(u,v)= U V$ where\\
 $U = \displaystyle\frac{1}{u_{1}}\left[ \frac{u_{1}u_{2}}{b_{1}b_{2} . . . b_{n}}(1,b_{0})\prod_{i=2}^{n}M(b_{i-1},x_{i},b_{i})\left( \begin{array}{cccccc}
1 \\
b_{n+1} \\
       
\end{array} \right)\right]^{\dfrac{1}{2}}$ and \\
$V = \displaystyle\frac{1}{u_{2}}\left[ \frac{u_{1}u_{2}}{b_{1}b_{2} . . . b_{n}}(1,b_{0})\prod_{i=2}^{n}M(b_{i-1},x_{i},b_{i})\left( \begin{array}{cccccc}
1 \\
b_{n+1} \\
       
\end{array} \right)\right]^{\dfrac{1}{2}}$  for the diagonal ray with origin $u_{1}u_{2}$ or\\
 $U = \displaystyle\frac{1}{u_{n}}\left[ \frac{u_{n}u_{n+1}}{b_{1}b_{2} . . . b_{n}}(1,b_{0})\prod_{i=2}^{n}M(b_{i-1},x_{i},b_{i})\left( \begin{array}{cccccc}
1 \\
b_{n+1} \\
       
\end{array} \right)\right]^{\dfrac{1}{2}}$ and \\
$V = \displaystyle\frac{1}{u_{n+1}}\left[ \frac{u_{n}u_{n+1}}{b_{1}b_{2} . . . b_{n}}(1,b_{0})\prod_{i=2}^{n}M(b_{i-1},x_{i},b_{i})\left( \begin{array}{cccccc}
1 \\
b_{n+1} \\
       
\end{array} \right)\right]^{\dfrac{1}{2}}$  for the diagonal ray with origin $u_{n}u_{n+1}$.  $\square$\\
\end{Cor}

We compute the transjective cluster variables by this method. We give later an example to illustrate the above results.

\subsection{Computation of nontransjective cluster variables of cluster algebra of type $\tilde{\mathbb{D}}$}

According to theorem $3.9$ the correspondence between a frieze associated with a seed $\mathcal{G}$ of type $\tilde{\mathbb{D}}_{n}$ and an $SL_{2}$-tiling associated with a particular seed of type $\tilde{\mathbb{A}}_{2n-1}$, allows us to regard transjective variables of the cluster algebra $\mathcal{A}(\mathcal{G})$ of type $\tilde{\mathbb{D}}_{n}$ as those of a particular type $\tilde{\mathbb{A}}_{2n-1}$. This suggests the idea of looking for nontransjective variables of the cluster algebra $\mathcal{A}(\mathcal{G})$ also among cluster variables of the same cluster algebra of type $\tilde{\mathbb{A}}_{2n-1}$. To this end we introduce a property of $SL_{2}$-tilings from [BR10- formula $(5)$].\\

\begin{prop}
Given three successive columns $C_{0}, C_{1}$ and $C_{2}$ of an $SL_{2}$-tiling $t$, there is a unique coefficient $\alpha \in \mathbb{K}$ such that $C_{0} + C_{2} = \alpha C_{1}$. $\square$
\end{prop}

\hspace{0.3 cm} We call $\alpha$ the \textit{linearization coefficient} of column $C_{1}$. There is a similar property and definition for three successive rows $L_{0},L_{1}$ and $L_{2}$.\\

\hspace{0.3 cm} Let $a_{1},...,a_{n}$ be elements of some ring $R$, the signed continuant polynomials (see [AR12-2.1]) are defined recursively as follows:\\ $$q_{n}(a_{1},..,a_{n})= q_{n-1}(a_{1},..,a_{n-1})a_{n} - q_{n-2}(a_{1},..,a_{n-2}), \quad n\geq 1,$$ setting $q_{-1} = 0$ and $q_{0} = 1$. 

\hspace{0.3 cm} It is known according to [AR12-theorem $4.4$] that cluster variables of a cluster algebra of type $\tilde{\mathbb{A}}$ either appear as elements of the corresponding $SL_{2}$-tiling, or as continuant polynomials of the linearization coefficients of the $SL_{2}$-tiling. We recall here the technique of computation of the nontransjective cluster variables of cluster algebra of type $\tilde{\mathbb{A}}$ according to [AR12]. We use the linearization coefficients of columns, the case of rows being analogous. We need a new concept of \textit{word} associated with a set of successive columns.\\

\begin{defi}
Let $S = \{ C_{i},...,C_{j} \}$ be a finite set of successive columns of the tiling $t$ associated with a boundary $f$. We call \textit{word associated with the set $S$}, the portion of the boundary $f$ between its intersections with the columns $C_{i},...,C_{j}$ augmented with one step to the left and one to the right. 
\end{defi}
 
 \begin{exem}
 Consider the scheme of example $2.4$ and let $S$ be the set of successive columns containing the variables from $c_{0}$ to $c_{1}$. The word associated with this set $S$ of columns is: $c_{-1}xc_{0}xc_{1}yc_{2}yc_{3}xc_{4}$.\\
\end{exem} 

\hspace{0.3 cm} The next theorem from [AR12-4.8] gives an expression for the signed continuant polynomial of the linearization coefficients of the tiling.\\
 
 \begin{theo}
 Consider an $SL_{2}$-tiling $t$ associated with some boundary $f$ with variables in $\mathbb{K}$.
 Let $C_{1},...,C_{k}$ be $k$ successive columns of the tiling $t$, with linearization coefficients $\alpha_{1},..., \alpha_{k}$. Let $m = b_{0}x_{1}b_{1} x_{2}...b_{n}x_{n+1}b_{n+1}$, $n \geq 1$, $x_{i} \in \{ x,y \}$, $ b_{i} \in \mathbb{K}$ be the word associated with this set of columns. Then the signed continuant polynomial $q_{k}(\alpha_{1},...,\alpha_{k})$ is equal to\\   $\displaystyle\frac{1}{b_{1}b_{2} . . . b_{n}}(b_{0}, 1)\prod_{i=2}^{n}M(b_{i-1},x_{i},b_{i})\left( \begin{array}{ccc}
1 \\
b_{n+1} \\
       
\end{array} \right)$. $\square$\\
 \end{theo}

\hspace{0.3 cm} We recall that in the case $\tilde{\mathbb{A}}_{r,s}$ with $r$ clockwise oriented arrows and $s$ anti-clockwise oriented arrows, the nontransjective cluster variables are lying on two tubes of rank $r$ and $s$ (see [SS-XIII.2.2]). In case $\tilde{\mathbb{D}}_{n}$ nontransjective variables lie on three tubes such that two are of rank $2$ and the third one is of rank $(n-2)$ (see [SS-XIII.2.2]).
 
\hspace{0.3 cm} Consider our quiver $Q'$ of type  $\tilde{\mathbb{A}}_{2n-1} = \tilde{\mathbb{A}}_{r,s}$, with $r+s = 2n$, which is of the form:\\ \resizebox{6cm}{!}{
$
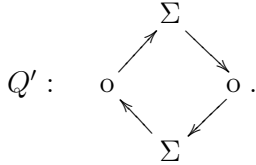
\xymatrix@=1pc{ 
 &   & &  \Sigma\ar[dr]  & &&&        \\
   &Q':   & \text{o}\ar[ur]   &      & \text{o} \; . \ar[dl]   & &&        \\
 &   &    &   \Sigma\ar[ul]&    &      &      &    \\
}$}\\
\hspace{0.3 cm} The analysis of orientation of arrows of $Q'$ gives us the equation $r = s+4$. From this and the equation $r+s = 2n$ we obtain $2s+4 = 2n$ and then $s = n-2$. We shall prove that the tube of rank $s$ in a particular case $\tilde{\mathbb{A}}_{2n-1}$ corresponding to the seed $\mathcal{G}'$ with $u_{0}=1$ can be identified with the tube of rank $(n-2)$ in the case $\tilde{\mathbb{D}}_{n}$.

\hspace{0.3 cm} According to [D10-theorem $5.1$], it is enough to know the variables lying on the mouth of a tube to determine all the variables lying on this tube. \\

\hspace{0.3 cm} We prove the following theorem which allows us to compute the $(n-2)$ nontransjective cluster variables lying on the mouth of the tube of rank $(n-2)$ in the case $\tilde{\mathbb{D}}_{n}$.\\

\begin{theo}
Let $n \geq 4$ be an integer, $\mathcal{G}$ a seed of type $\tilde{\mathbb{D}}_{n}$ and $\mathcal{G}'$ the seed of type $\tilde{\mathbb{A}}_{2n-1}$ associated with $\mathcal{G}$.
The tube of rank $s=(n-2)$ containing nontransjective variables of the cluster algebra $\mathcal{A}(\mathcal{G}')$ with $u_{0}=1$ coincides with the tube of rank $(n-2)$ containing nontransjective cluster variables of the cluster algebra $\mathcal{A}(\mathcal{G})$ of type $\tilde{\mathbb{D}}_{n}$. 
\end{theo}

\hspace{0.3 cm} Before the proof of theorem $4.6$, we recall formula $(1.1)$ from [ADSS12-1.3], which we are going to use to compute variables lying on the mouths of tubes. For any locally finite quiver $Q$ with variables (here $Q$ is a quiver of type $\tilde{\mathbb{D}}_{n}$), we define a family of matrices with coefficients in $\mathbb{Z}[u_{i} | i \in Q_{0}]$ as follows:\\
for any arrow $s(\epsilon)\stackrel{\epsilon}{\rightarrow} t(\epsilon) $ in $Q_{1}$ we set:\\

$ M(\epsilon)= \left( \begin{array}{cccccc}
u_{t(\epsilon)}&0  \\
1 &u_{s(\epsilon)} \\
\end{array} \right)\quad$ and  $\quad  M(\epsilon^{-1})= \left( \begin{array}{cccccc}
u_{t(\epsilon)} &1   \\
0  &u_{s(\epsilon)} \\
\end{array} \right)$. \\

Consider now, for all $k \in \{0,1,...n \}$, a reduced walk $c=v_{1}\stackrel{d_{1}}{-} \ldots \stackrel{d_{m}}{-} v_{m+1} $ from a vertex $v_{1}$ to a vertex $v_{m+1}$ in $Q$ (for the notion of reduced walk, we refer to [ASS-II.1.1]), here $d_{l}$, are arrows.\\

 We define $$V_{c}(k)= \left[ \begin{array}{cccccc}\displaystyle
\prod u_{t(\epsilon)} \atop {\epsilon \in Q_1(v_{k},-)\atop \epsilon \neq d^{\small{\pm 1}}_{k},d^{\small{\pm 1}}_{k-1}} &0  \\
0 &\displaystyle \prod u_{s(\epsilon)} \atop {\epsilon \in Q_1(-,v_{k})\atop \epsilon \neq d^{\small{\pm 1}}_{k},d^{\small{\pm 1}}_{k-1}} \\
\end{array}\right],$$ where $Q_{1}(v_{k},-) = \{ \epsilon \in Q_{1}| s(\epsilon)=v_{k} \}, \, \, Q_{1}(-, v_{k}) = \{ \epsilon \in Q_{1}| t(\epsilon)=v_{k} \}$. The empty product is equal to $1$.

\hspace{0.3 cm} There is [ADSS12-4.1] a unique cluster variable of $\mathcal{A}(\mathcal{G})$ corresponding to every reduced walk $c$ of length $m$ in $Q$. This cluster variable is given by ([ADSS12-1.3]):

$$\qquad \qquad \qquad \displaystyle\frac{1}{\prod_{k=0}^{m} u_{t(d_{k})}}\left[ 1, 1\right] \left( \prod_{k=0}^{m}M(d_{k})V_{c}(k+1)\right) \left[  \begin{array}{cccccc}
1 \\
1 \\
       
\end{array} \right] \qquad \qquad \qquad  (6)$$ where $M(d_{0})$ is the identity matrix by convention. \\

\textbf{Proof of theorem 4.6}

We give the proof for the quiver $Q$ with the following orientation:

$$
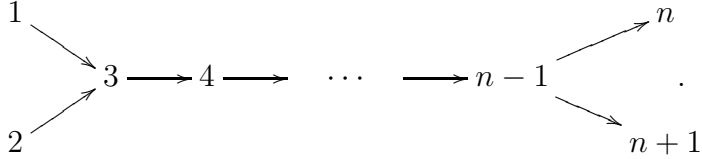
\xymatrix@R10pt{
1\ar[rd]&&&&&n\\
&3\ar[r]&4\ar[r]&\quad\cdots\quad\ar[r]&n-1\ar[ru]\ar[rd]& \quad .\\
2\ar[ru]&&&&&n+1
}$$ 

\hspace{0.3 cm} Due to the fact that all other orientations of the $\tilde{\mathbb{D}}_{n}$ diagram are mutation-equivalent to the chosen one (because of lemma $3.1$), the proof for other orientations can be reduced to this case.
 
\hspace{0.3 cm} It is enough to prove, according to [D10-5.1], that the variables lying on the mouth of the tube of rank $s=(n-2)$ in case $\tilde{\mathbb{A}}_{2n-1}$ corresponding to the seed $\mathcal{G}'$ with $u_{0}=1$ are all the $(n-2)$ nontransjective cluster variables lying on the mouth of the tube of rank $(n-2)$ in the case of the seed $\mathcal{G}$ of type $\tilde{\mathbb{D}}_{n}$.

\hspace{0.3 cm} According to [AD11-1.3], the nontransjective cluster variables  of cluster algebra of type $\tilde{\mathbb{D}}_{n}$ lying on the mouth of the tube of rank $(n-2)$ correspond to simple $\mathbb{K}Q$-modules $S_{i}$ associated with points between the two joints of forks inclusively, that is with the points $i = 3,...,n-1$. Then we obtain $(n-3)$ simple $\mathbb{K}Q$-modules lying on the mouth of the tube.\\

\hspace{0.3 cm} To compute the cluster variables of the cluster algebra $\mathcal{A}(\mathcal{G})$ of type $\tilde{\mathbb{D}}_{n}$ corresponding to these simple modules, we identify each point with a walk of length $0$ and apply formula $(6)$. It is easy to see that the obtained variables are equal to the linearization coefficients of columns of the tiling $t$ passing through the vertices $u_{3},...,u_{n-1}$ of the root $\bar{\omega}$ for our orientation of $Q$. On the other hand we know [AR12] that the linearization coefficients of rows and columns of the $SL_{2}$-tiling $t$ below the boundary $\tilde{f}_{0}$ are the variables lying on the mouths of tubes of rank $r$ and $s$, respectively, of the Auslander-Reiten quiver corresponding to $Q'$ of type $\tilde{\mathbb{A}}_{2n-1}$. Therefore we conclude that the $(n-3)$ cluster variables of cluster algebra of type $\tilde{\mathbb{D}}_{n}$ obtained on the mouth of the tube of rank $(n-2)$ coincide with the $(n-3)$ variables on the mouth of the tube of rank $s=(n-2)$ in case $\tilde{\mathbb{A}}_{2n-1}$ and are set in the same order.\\

\hspace{0.3 cm} Thus we have found $(n-3)$ out of $(n-2)$ cluster variables lying on the mouth of the tube in the case  $\mathcal{A}(\mathcal{G})$ of type $\tilde{\mathbb{D}}_{n}$. To compute the $(n-2)$th cluster variable we use the existence of linear recurrence relations between the corresponding values of the frieze $F$ at a joint of a fork and its neighbours from [KS11-6]. 

\hspace{0.3 cm} Namely, we denote by $\tau$ the Auslander-Reiten translation (see [ASS-IV.2.3]), $X^{i}_{k} = \mathfrak{a}(k,i)$ and by $X_{\tau^{k}S_{i}}$ (note that all cluster variables lying on the mouth of a tube of rank $(n-2)$ are of this form) the image by the Caldero-Chapoton map (see [CC06]) of the $\mathbb{K}Q$-module $\tau^{k}S_{i}$. Corresponding to the sequence $X_{k}^{3}$ associated with the joint $3$ of the quiver $Q$, we get the following relation: $$ X^{4}_{k} = X_{k}^{3}X_{\tau^{k}S_{4}} - X_{k}^{2}X^{1}_{k}. $$

\hspace{0.3 cm} Rewriting this relation as follows: $  X_{k}^{3}X_{\tau^{k}S_{4}} = X_{k}^{4} + X_{k}^{2}X^{1}_{k} $, one can see that each variable $X_{\tau^{k}S_{4}}$ lying on the mouth of the tube is the linearization coefficient of a column passing through the variable $X_{k}^{3}$ in the $SL_{2}$-tiling $t$ defined by the boundary $\tilde{f}_{0}$. This allows us to conclude that the $(n-2)${th} cluster variable also corresponds to a linearization coefficient in the $SL_{2}$-tiling $t$.

\hspace{0.3 cm} Thus the $(n-2)$ cluster variables of the cluster algebra $\mathcal{A}(\mathcal{G})$ of type $\tilde{\mathbb{D}}_{n}$ lying on the mouth of the tube of rank $(n-2)$ are equal to the linearization coefficients of columns passing through the points in the tiling below the boundary $\tilde{f}_{0}$ which correspond to the joints of the forks in the quiver ${F}$. (Recall [AR12] that the sequence of linearization coefficients of the columns of the $SL_{2}$-tiling corresponding to a quiver of type $\tilde{\mathbb{A}}_{r,s}$ is periodic with the period $s$. Since in our case $s = n-2$, there are only $(n-2)$ such coefficients.) 

\hspace{0.3 cm} Thus we identify the tube of rank $s=(n-2)$ in the case $\tilde{\mathbb{A}}_{2n-1}$ (corresponding to the cluster algebra $\mathcal{A}(\mathcal{G}')$ with $u_{0}=1$) with the tube of rank $(n-2)$ in case $\tilde{\mathbb{D}}_{n}$ (corresponding to the cluster algebra $\mathcal{A}(\mathcal{G})$). $\square$ \\

\begin{rem}

Let $q_{n}$ be the signed continuant polynomial, ${a}_{k,n}$ with $k, n \, \in \mathbb{N}$, be a variable of depth $n$ on the tube and ${a}_{k,1}$ be a variable lying on the mouth of the tube. According to [AD11-3.5] the variables lying on a tube with rank $p \geq 1$ are related by the following relation ${a}_{i,n} = q_{n}\left( {a}_{i,1},...,{a}_{i+n-1,1}\right) $.\\

\end{rem}

\hspace{0.3 cm} It remains now to determine the nontransjective cluster variables lying on the mouths of the two tubes of rank $2$.\\

\hspace{0.3 cm} We do this by interpreting the results on quiver representations from [SS-XIII-2.6-a and b] in terms of cluster variables. This gives that the nontransjective cluster variables on the mouth of each tube of rank $2$ correspond to reduced walks in the seed $\mathcal{G}$ as follows:\\
\begin{enumerate}
\item For the first tube
\begin{enumerate}
\item the unique reduced walk from $u_{1}$ to $u_{n+1}$
\item the unique reduced walk from $u_{2}$ to $u_{n}$
\end{enumerate}
\item For the second tube
\begin{enumerate}
\item the unique reduced walk from $u_{1}$ to $u_{n}$
\item the unique reduced walk from $u_{2}$ to $u_{n+1}$
\end{enumerate}
\end{enumerate}  

\hspace{0.3 cm} Then the nontransjective cluster variables of cluster algebra of type $\tilde{\mathbb{D}}_{n}$ lying on the tubes of rank $2$ are obtained by applying to these walks formula $(6)$.

\hspace{0.3 cm} We are now able to compute the transjective and nontransjective cluster variables of cluster algebra of type $\tilde{\mathbb{D}}_{n}$.\\

\hspace{0.3 cm} We give now an example to illustrate the above results.

\begin{exem}
Consider the quiver $Q$ of type $\tilde{\mathbb{D}}_{4}$ of example $3.6$.\\
$\xymatrix @=10pt{
& \text{We have} \; \; \textbf{$\Sigma$} :\;{1}\ar[r]
&{3}
&4\ar[l] 
}  $
$\xymatrix @=10pt{
& \text{and} \; \; \textbf{$^{t}\Sigma$} :\;{4}\ar[r]
&{3}
&1\ar[l] 
}$.\\ The root $\bar{\omega}$ of $Q$ with variables is: $\bar{\omega} = u_{1}u_{2}xu_{3}yu_{4}u_{5}$. The boundary $\tilde{f}_{0}$ associated with $Q$ is:
$\tilde{f}_{0} = ^{\infty}(u_{1}u_{2}xu_{3}yu_{4}u_{5}x1xu_{4}u_{5}xu_{3}yu_{1}u_{2}x1x)(u_{1}u_{2}xu_{3}yu_{4}u_{5})y1\\yu_{4}u_{5}xu_{3}yu_{1}u_{2}y1yu_{1}u_{2}xu_{3}yu_{4}u_{5})^{\infty}$, this gives the following $SL_{2}$-tiling below the boundary $\tilde{f}_{0}$: \\

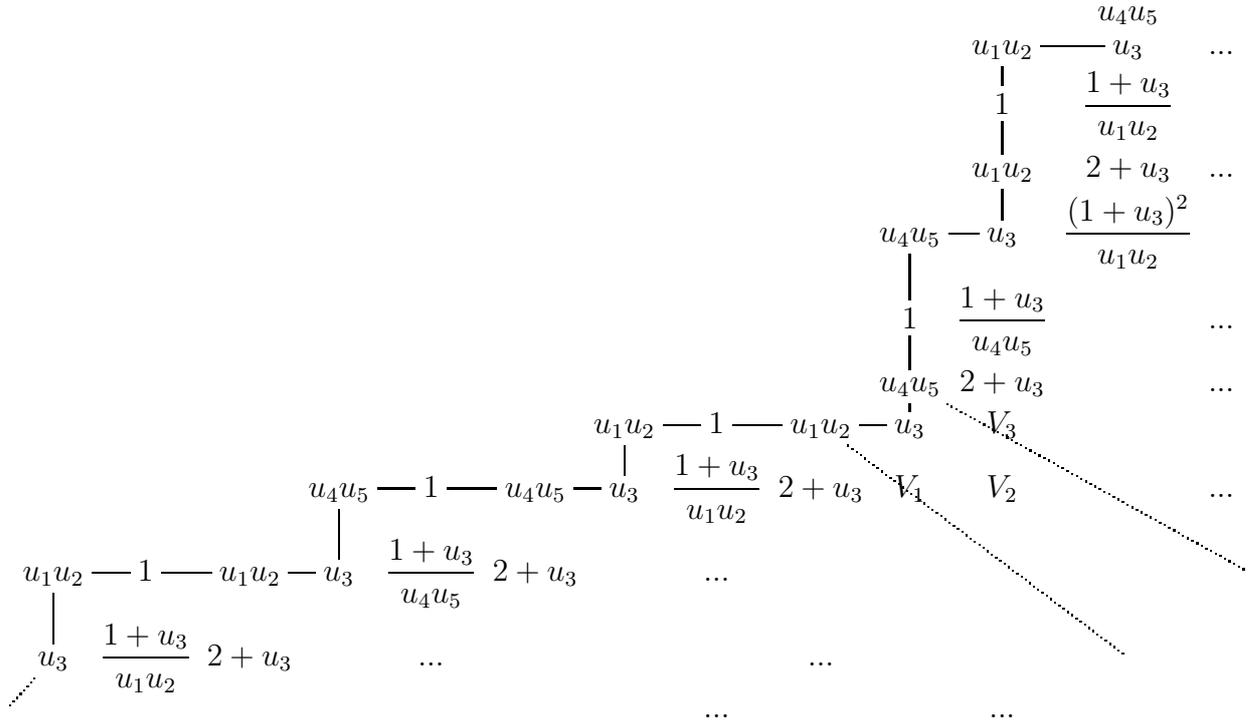
$$\xymatrix@=0pt{
&&&&&&&&&&&&u_{4}u_{5}&&\\
&&&&&&&&&&&u_{1}u_{2}\ar@{-}[r]&u_{3}\ar@{-}[u]&...&&\\
&&&&&&&&&&&1\ar@{-}[u]&\displaystyle\frac{1+u_{3}}{u_{1}u_{2}}&&&&&\\
&&&&&&&&&&&u_{1}u_{2}\ar@{-}[u]&2+u_{3}&...&&&&\\
&&&&&& &   &   &    &  u_{4}u_{5}\ar@{-}[r] &   u_{3} \ar@{-}[u]&  \displaystyle\frac{(1+u_{3})^{2}}{u_{1}u_{2}}   &  &               \\
&&&&& &   &   &  &  &   1 \ar@{-}[u]&  \displaystyle\frac{1+u_{3}}{u_{4}u_{5}}&   &...&       \\
&&&&& &&   &  &  &u_{4}u_{5}\ar@{-}[u]\ar@{.}[dddrrrr]& 2+u_{3}&  &... &   &          \\
&&&&&&&u_{1}u_{2}\ar@{-}[r]&1\ar@{-}[r] &u_{1}u_{2}\ar@{-}[r]\ar@{.}[dddrrr] & u_{3}\ar@{-}[u]&V_{3}& &&  &&&&      \\
&&&&u_{4}u_{5}\ar@{-}[r]&1\ar@{-}[r]&u_{4}u_{5}\ar@{-}[r]&u_{3}\ar@{-}[u]&\displaystyle\frac{1+u_{3}}{u_{1}u_{2}}&2+u_{3}&V_{1}&V_{2}&&... & &         \\
&u_{1}u_{2}\ar@{-}[r]&1\ar@{-}[r]&u_{1}u_{2}\ar@{-}[r]&u_{3}\ar@{-}[u] &\displaystyle\frac{1+u_{3}}{u_{4}u_{5}}&2+u_{3}& & ...& & && &&          \\
&u_{3}\ar@{.}[dl]\ar@{-}[u]&\displaystyle\frac{1+u_{3}}{u_{1}u_{2}}&2+u_{3}&&...&&&&...&&&&&&&&&\\
&&&&& & &&...&&&...&\\
}$$\\

\hspace{0.3 cm} The oblique band delimited by the diagonal rays with origins $u_{1}u_{2}$ and $u_{4}u_{5}$ corresponds, according to theorem $3.9$, to the modelled quiver $\bar{F}$ and contains the transjective cluster variables of cluster algebra of type $\tilde{\mathbb{D}}_{4}$.

\hspace{0.3 cm} The variables $V_{1}$ and $V_{3}$ represent the products of two cluster variables on the bottom and upper lines in the modelled quiver $\bar{F}$ and $V_{2}$ represents some cluster variable in $\bar{F}$.\\

\hspace{0.3 cm} The cluster variables in these three positions are computed as follows.\\

\hspace{0.3 cm} The words associated with $V_{1}$, $V_{2}$ and $V_{3}$ are $u_{3}yu_{1}u_{2}x1xu_{1}u_{2}xu_{3}$,\\ $u_{3}yu_{1}u_{2}x1xu_{1}u_{2}xu_{3}yu_{4}u_{5}y1yu_{4}u_{5}xu_{3}$ and $u_{3}yu_{4}u_{5}y1yu_{4}u_{5}xu_{3}$.\\

Then we have, by applying theorem $2.5$, the following results:\\
\begin{eqnarray}
V_{1} &=& \displaystyle\frac{1}{u_{1}^{2}u_{2}^{2}} \left( 1,u_{3}\right)  \left( \begin{array}{cccccc}
u_{1}u_{2} &  &1 \\
0 &  &1 \\
\end{array} \right) \left( \begin{array}{cccccc}
1 &  &1 \\
0 &  &u_{1}u_{2} \\
\end{array} \right) \left( \begin{array}{cccccc}
1 \\
u_{3} \\     
\end{array} \right) \nonumber \\
&=&\displaystyle\frac{(1+u_{3})^{2}}{u_{1}u_{2}} \nonumber\\
&=&\displaystyle\frac{1+u_{3}}{u_{1}}\times \frac{1+u_{3}}{u_{2}}.\nonumber
\end{eqnarray}

\begin{eqnarray}
V_{3} &=& \displaystyle\frac{1}{u_{4}^{2}u_{5}^{2}} \left( 1,u_{3}\right)  \left( \begin{array}{cccccc}
1 &  &0 \\
1 &  &u_{4}u_{5} \\
\end{array} \right) \left( \begin{array}{cccccc}
u_{4}u_{5} &  &0 \\
1 &  &1 \\
\end{array} \right) \left( \begin{array}{cccccc}
1 \\
u_{3} \\     
\end{array} \right) \nonumber \\
&=&\displaystyle\frac{(1+u_{3})^{2}}{u_{4}u_{5}} \nonumber\\
&=&\displaystyle\frac{1+u_{3}}{u_{4}}\times \frac{1+u_{3}}{u_{5}}.\nonumber
\end{eqnarray}

\hspace{0.3 cm} $V_{1}$ and $V_{3}$ being placed on the bottom and upper lines in $\bar{F}$, we write these values in the form of a product of two transjective cluster variables as in corollary $4.1$ above.  The other position corresponds to the following transjective cluster variable of type $\tilde{\mathbb{D}}_{4}$:\\ 

\begin{eqnarray}
V_{2} &=& \displaystyle\frac{1}{(u_{1}u_{2})^{2}u_{3}(u_{4}u_{5})^{2}} \left( 1,u_{3}\right)  \left( \begin{array}{cccccc}
u_{1}u_{2}  &1  \\
0 & 1 \\
\end{array} \right) \left( \begin{array}{cccccc}
1 &1  \\
0 & u_{1}u_{2} \\
\end{array} \right) \left( \begin{array}{cccccc}
u_{1}u_{2} & 1   \\
0 & u_{3} \\
\end{array} \right) \times \nonumber \\ & \times & \left( \begin{array}{cccccc}
u_{4}u_{5} &0   \\
1& u_{3} \\
\end{array} \right)  \left( \begin{array}{cccccc}
1 & 0   \\
1 & u_{4}u_{5} \\
\end{array} \right) \left( \begin{array}{cccccc}
u_{4}u_{5} & 0   \\
1 & 1 \\
\end{array} \right) \left( \begin{array}{cccccc}
1 \\
u_{3} \\     
\end{array} \right) \nonumber \\
 &=& \displaystyle\frac{1}{(u_{1}u_{2})^{2}u_{3}(u_{4}u_{5})^{2}} \left( 1,u_{3}\right)  \left( \begin{array}{cccccc}
(u_{1}u_{2})^{2}  &u_{1}u_{2}(1+2u_{3})  \\
0 &  u_{1}u_{2}u_{3}\\
\end{array} \right)  \times \nonumber \\ & \times & \left( \begin{array}{cccccc}
(u_{4}u_{5})^{2} & 0   \\
u_{4}u_{5}(1+2u_{3}) & u_{3}u_{4}u_{5} \\
\end{array} \right) \left( \begin{array}{cccccc}
1 \\
u_{3} \\     
\end{array} \right) \nonumber \\
 &=& \displaystyle\frac{1}{(u_{1}u_{2})^{2}u_{3}(u_{4}u_{5})^{2}} \left( (u_{1}u_{2})^{2}  ,u_{1}u_{2}(1+u_{3})^{2}\right)   \left( \begin{array}{cccccc}
(u_{4}u_{5})^{2} \\
u_{4}u_{5}(1+u_{3})^{2} \\     
\end{array} \right)  \nonumber \\
&=&\displaystyle\frac{u_{1}u_{2}u_{3}u_{4}u_{5}+(1+u_{3})^{4}}{u_{1}u_{2}u_{3}u_{4}u_{5}}.\nonumber
\end{eqnarray}
 
\hspace{0.3 cm} We compute thus all transjective cluster variables of cluster algebra of type $\tilde{\mathbb{D}}_{4}$ by this method.\\
 
\hspace{0.3 cm} We compute now nontransjective cluster variables of cluster algebra of type $\tilde{\mathbb{D}}_{4}$. Because $n = 4$, these cluster variables all lie in tubes of rank $2$.

\hspace{0.3 cm} For one of these tubes, say $\mathcal{T}_{1}$, we have two nontransjective cluster variables $\alpha_{1}$ and $\alpha'_{1}$ on the mouth. These variables are given by the linearization coefficients of columns of the $SL_{2}$-tiling passing through points corresponding to the joints of forks in $F$. These cluster variables can also be obtained by applying the formula of theorem $4.5$. 

\hspace{0.3 cm} We choose to use theorem $4.5$ to compute the linearization coefficient $\alpha_{1}$ of the column passing through the vertex $u_{3}$ of the root $\bar{\omega}$ and to use proposition $4.2$ to compute the linearization coefficient $\alpha'_{1}$ of the column passing through $V_{2}$.\\

\hspace{0.3 cm} For the computation of $\alpha_{1}$ we apply theorem $4.5$. That is, first we need to determine the word associated with the column containing the vertex $u_{3}$ of the root $\bar{\omega}$. According to the embedding of the boundary $\tilde{f}_{0}$ in the plane, this word is:\\ $u_{1}u_{2}xu_{3}yu_{4}u_{5}y1yu_{4}u_{5}xu_{3}.$ Applying the formula of theorem $4.5$ to this word we obtain:\\

\begin{eqnarray}
\alpha_{1} &=& \displaystyle\frac{1}{u_{3}(u_{4}u_{5})^{2}} \left( u_{1}u_{2},1 \right)   \left( \begin{array}{cccccc}
u_{4}u_{5} & 0   \\
1&  u_{3}\\
\end{array} \right)  \left( \begin{array}{cccccc}
1 &  0   \\
1 &  u_{4}u_{5}\\
\end{array} \right) \left( \begin{array}{cccccc}
u_{4}u_{5} &0   \\
1 &1 \\
\end{array} \right) \left( \begin{array}{cccccc}
1 \\
u_{3} \\     
\end{array} \right) \nonumber \\
 &=& \displaystyle\frac{1}{u_{3}(u_{4}u_{5})^{2}} \left( 1+u_{1}u_{2}u_{4}u_{5} , u_{3}\right)  \left( \begin{array}{cccccc}
1 &  0   \\
1 &  u_{4}u_{5} \\
\end{array} \right)  \left( \begin{array}{cccccc}
u_{4}u_{5} \\
1+u_{3} \\     
\end{array} \right) \nonumber \\
 &=& \displaystyle\frac{1}{u_{3}(u_{4}u_{5})^{2}} \left( 1+u_{1}u_{2}u_{4}u_{5}  ,u_{3}\right)   \left( \begin{array}{cccccc}
u_{4}u_{5} \\
u_{4}u_{5}(2+u_{3}) \\     
\end{array} \right)  \nonumber \\
&=&\displaystyle\frac{u_{1}u_{2}u_{4}u_{5}+(1+u_{3})^{2}}{u_{3}u_{4}u_{5}}.\nonumber
\end{eqnarray}

\hspace{0.3 cm} We apply proposition $4.2$ to compute the linearization coefficient $\alpha'_{1}$ of the column passing through the variable $V_{2}$. Let us calculate $\alpha'_{1}$ at the vertex $u_{3}$ at the intersection of the column containing $V_{2}$ with the boundary $\tilde{f}_{0}$. Then we have:\\
 
$\alpha'_{1} = \displaystyle\frac{1}{u_{3}}\left[ \frac{(1+u_{3})^{2}}{u_{1}u_{2}} + u_{4}u_{5} \right] = \frac{ u_{1}u_{2}u_{4}u_{5} + (1+u_{3})^{2} }{u_{1}u_{2}u_{3}} $.\\

\hspace{0.3 cm} Thus, $\alpha_{1}$ and $\alpha'_{1}$ are nontransjective cluster variables lying on the mouth of the tube $\mathcal{T}_{1}$.\\

\hspace{0.3 cm} For the second tube $\mathcal{T}_{2}$, we associate the modules lying on the mouth with the following reduced walks from $u_{1}$ to $u_{5}$ and from $u_{2}$ to $u_{4}$: $u_{1}xu_{3}yu_{5}$ and $u_{2}xu_{3}yu_{4}$. By applying formula $(6)$, we obtain the corresponding nontransjective cluster variables $\alpha_{2}$ and $\alpha'_{2}$:

\begin{eqnarray}
\alpha_{2} &=& \displaystyle\frac{1}{u_{1}u_{3}u_{5}} \left( 1, 1 \right)  \left( \begin{array}{cccccc}
1 &0  \\
0 &1 \\
\end{array} \right)  \left( \begin{array}{cccccc}
u_{3} & 0   \\
1&  u_{1}\\
\end{array} \right)  \left( \begin{array}{cccccc}
1 &  0   \\
0 &  u_{2}u_{4}\\
\end{array} \right) \times \nonumber \\
&\times & \left( \begin{array}{cccccc}
u_{3} &1   \\
0 &u_{5} \\
\end{array} \right)\left( \begin{array}{cccccc}
1 &0   \\
0 &1 \\
\end{array} \right)  \left( \begin{array}{cccccc}
1 \\
1 \\     
\end{array} \right) \nonumber \\
 &=& \displaystyle\frac{1}{u_{1}u_{3}u_{5}} \left( 1+u_{3} , u_{1}\right)  \left( \begin{array}{cccccc}
1 &  0   \\
0 &  u_{2}u_{4} \\
\end{array} \right)  \left( \begin{array}{cccccc}
1 + u_{3} \\
u_{5} \\     
\end{array} \right) \nonumber \\
 &=& \displaystyle\frac{1}{u_{1}u_{3}u_{5}} \left( 1 + u_{3}  ,u_{1}\right)   \left( \begin{array}{cccccc}
1 + u_{3} \\
u_{2}u_{4}u_{5} \\     
\end{array} \right)  \nonumber \\
&=&\displaystyle\frac{u_{1}u_{2}u_{4}u_{5}+(1+u_{3})^{2}}{u_{1}u_{3}u_{5}}.\nonumber
\end{eqnarray}

\begin{eqnarray}
\alpha'_{2} &=& \displaystyle\frac{1}{u_{2}u_{3}u_{4}} \left( 1, 1 \right)  \left( \begin{array}{cccccc}
1 &0  \\
0 &1 \\
\end{array} \right)  \left( \begin{array}{cccccc}
u_{3} & 0   \\
1&  u_{2}\\
\end{array} \right)  \left( \begin{array}{cccccc}
1 &  0   \\
0 &  u_{1}u_{5}\\
\end{array} \right) \times \nonumber \\
&\times & \left( \begin{array}{cccccc}
u_{3} &1   \\
0 &u_{4} \\
\end{array} \right)\left( \begin{array}{cccccc}
1 &0   \\
0 &1 \\
\end{array} \right)  \left( \begin{array}{cccccc}
1 \\
1 \\     
\end{array} \right) \nonumber \\
 &=& \displaystyle\frac{1}{u_{2}u_{3}u_{4}} \left( 1+u_{3} , u_{2}\right)  \left( \begin{array}{cccccc}
1 &  0   \\
0 &  u_{1}u_{5} \\
\end{array} \right)  \left( \begin{array}{cccccc}
1 + u_{3} \\
u_{4} \\     
\end{array} \right) \nonumber \\
 &=& \displaystyle\frac{1}{u_{2}u_{3}u_{4}} \left( 1 + u_{3}  ,u_{2}\right)   \left( \begin{array}{cccccc}
1 + u_{3} \\
u_{1}u_{4}u_{5} \\     
\end{array} \right)  \nonumber \\
&=&\displaystyle\frac{u_{1}u_{2}u_{4}u_{5}+(1+u_{3})^{2}}{u_{2}u_{3}u_{4}}.\nonumber
\end{eqnarray}

\hspace{0.3 cm} Thus $\alpha_{2}$ and $\alpha'_{2}$ are the nontransjective cluster variables lying on the mouth of $\mathcal{T}_{2}$.\\
   
\hspace{0.3 cm} For the third tube $\mathcal{T}_{3}$, we associate the modules lying on the mouth with the following reduced walks from $u_{1}$ to $u_{4}$ and from $u_{2}$ to $u_{5}$: $u_{1}xu_{3}yu_{4}$ and $u_{2}xu_{3}yu_{5}$. By applying formula $(6)$, we obtain the corresponding nontransjective cluster variables $\alpha_{3}$ and $\alpha'_{3}$:

\begin{eqnarray}
\alpha_{3} &=& \displaystyle\frac{1}{u_{1}u_{3}u_{4}} \left( 1, 1 \right)  \left( \begin{array}{cccccc}
1 &0  \\
0 &1 \\
\end{array} \right)  \left( \begin{array}{cccccc}
u_{3} & 0   \\
1&  u_{1}\\
\end{array} \right)  \left( \begin{array}{cccccc}
1 &  0   \\
0 &  u_{2}u_{5}\\
\end{array} \right) \times \nonumber \\
&\times & \left( \begin{array}{cccccc}
u_{3} &1   \\
0 &u_{4} \\
\end{array} \right)\left( \begin{array}{cccccc}
1 &0   \\
0 &1 \\
\end{array} \right)  \left( \begin{array}{cccccc}
1 \\
1 \\     
\end{array} \right) \nonumber \\
 &=& \displaystyle\frac{1}{u_{1}u_{3}u_{4}} \left( 1+u_{3} , u_{1}\right)  \left( \begin{array}{cccccc}
1 &  0   \\
0 &  u_{2}u_{5} \\
\end{array} \right)  \left( \begin{array}{cccccc}
1 + u_{3} \\
u_{4} \\     
\end{array} \right) \nonumber \\
 &=& \displaystyle\frac{1}{u_{1}u_{3}u_{4}} \left( 1 + u_{3}  ,u_{1}\right)   \left( \begin{array}{cccccc}
1 + u_{3} \\
u_{2}u_{4}u_{5} \\     
\end{array} \right)  \nonumber \\
&=&\displaystyle\frac{u_{1}u_{2}u_{4}u_{5}+(1+u_{3})^{2}}{u_{1}u_{3}u_{4}}.\nonumber
\end{eqnarray}

\begin{eqnarray}
\alpha'_{3} &=& \displaystyle\frac{1}{u_{2}u_{3}u_{5}} \left( 1, 1 \right)  \left( \begin{array}{cccccc}
1 &0  \\
0 &1 \\
\end{array} \right)  \left( \begin{array}{cccccc}
u_{3} & 0   \\
1&  u_{2}\\
\end{array} \right)  \left( \begin{array}{cccccc}
1 &  0   \\
0 &  u_{1}u_{4}\\
\end{array} \right) \times \nonumber \\
&\times & \left( \begin{array}{cccccc}
u_{3} &1   \\
0 &u_{5} \\
\end{array} \right)\left( \begin{array}{cccccc}
1 &0   \\
0 &1 \\
\end{array} \right)  \left( \begin{array}{cccccc}
1 \\
1 \\     
\end{array} \right) \nonumber \\
 &=& \displaystyle\frac{1}{u_{2}u_{3}u_{5}} \left( 1+u_{3} , u_{2}\right)  \left( \begin{array}{cccccc}
1 &  0   \\
0 &  u_{1}u_{4} \\
\end{array} \right)  \left( \begin{array}{cccccc}
1 + u_{3} \\
u_{5} \\     
\end{array} \right) \nonumber \\
 &=& \displaystyle\frac{1}{u_{2}u_{3}u_{5}} \left( 1 + u_{3}  ,u_{2}\right)   \left( \begin{array}{cccccc}
1 + u_{3} \\
u_{1}u_{4}u_{5} \\     
\end{array} \right)  \nonumber \\
&=&\displaystyle\frac{u_{1}u_{2}u_{4}u_{5}+(1+u_{3})^{2}}{u_{2}u_{3}u_{5}}.\nonumber
\end{eqnarray}

\hspace{0.3 cm} Thus $\alpha_{3}$ and $\alpha'_{3}$ are the nontransjective cluster variables lying on the mouth of the tube $\mathcal{T}_{3}$.\\ 
\end{exem}

ACKNOWLEDGEMENTS. The author is grateful to Ibrahim Assem and Vasilisa Shramchenko for useful discussions and careful reading of the manuscript and also to the department of mathematics of the University of Sherbrooke where this work was done. The author gratefully acknowledges support from CIDA (Canadian International Development Agency) fellowship.

{Kodjo Essonana MAGNANI\\
   D\'epartement de math\'ematiques\\
   Universit\'e de Sherbrooke\\
2500, boul. de l'Universit\'e,\\
Sherbrooke, Qu\'ebec,
J1K 2R1\\
   Canada\\
   Kodjo.essonana.magnani@USherbrooke.ca}


\begin{thebibliography}{99}


\bibitem[{AD11}] {} {I. Assem and G. Dupont}, Friezes and the construction of the euclidean cluster variables, \emph{J. Pure and Applied Algebra} {\bf 215} (2011), 2322-2340.

\bibitem[ADSS12]{}  I. Assem, G. Dupont, R. Schiffler and D. Smith, Friezes, strings and
cluster variables, \emph{Glasgow J. Math.} {\bf 54} (2012), no. 1, 27-60.

\bibitem[AR12] {} I. Assem and C. Reutenauer, Mutating Seeds: types $A$ and $\tilde{A}$ \emph{Ann. Math. Blaise Pascal} {\bf 19} (2012), no. 1, 29-73.

\bibitem[ARS10] {} I. Assem, C. Reutenauer, and D. Smith, {Friezes}, \emph{Adv. in Math.} {\bf 225} (2010), 3134-3165.

\bibitem[{ASS}]{} I. Assem, D. Simson and A. Skowro\'nski, {\em Elements of
  the Representation Theory of Associative Algebras 1:
  Techniques of Representation Theory\/}, London Mathematical Society
Student Texts 65, Cambridge University Press, (2006).

\bibitem [BM12]{} K. Baur, R. J. Marsh, Categorification of a frieze pattern determinant, \emph{J. Combin. Theory Ser.} A 119 (2012), 1110-1122.

\bibitem [BMR08]{} A. Buan, R. Marsh and I. Reiten, Cluster mutation via quiver representations, \emph{Comment. Math. Helv.} 83 (2008), no. 1, 143-177.

\bibitem [BMR09]{} A. Buan, R. Marsh and I. Reiten, Denominators of cluster variables, \emph{ J. Lond. Math. Soc.} (2) 79 (2009), no. 3, 589-611.

\bibitem [BMRRT06]{} A. Buan, R. Marsh, M. Reineke, I. Reiten and G. Todorov, Tilting theory and cluster combinatorics, \emph{Adv. Math.} 204 (2006), no. 2, 572-618.

\bibitem [BR10]{} F. Bergeron, C. Reutenauer, $SL_{k}$-tilings of the plane, \emph{Illinois J. Math.} 54 (2010), no. 1, 263-300.

\bibitem [CC06]{} P. Caldero, F. Chapoton, Cluster algebras as Hall algebras of quiver representations, \emph{Commentarii Mathematici Helvetici} 81 (2006), 596-616.

\bibitem[{CCS06}]{}  { P. Caldero, F. Chapoton and
R. Schiffler}, Quivers with relations arising from clusters ($A_n$
case), \emph{Trans. Amer. Math. Soc.} {\bf 358} (2006), no. 3, 1347-1364. 

\bibitem [CC73-I]{} J. Conway, H. Coxeter, Triangulated polygons and frieze patterns, \emph{The Mathematical Gazette} 57 (400) (1973), 87-94. 

 \bibitem [CC73-II]{} J. Conway, H. Coxeter, Triangulated polygons and frieze patterns, \emph{The Mathematical Gazette} 57 (401) (1973), 175-183.
 
 \bibitem [C71]{} H. Coxeter, Frieze Patterns, \emph{Acta Arithmetica} XVIII (1971), 297-310.

\bibitem[{D10}]{} G. Dupont, Cluster multiplication in regular components via generalized Chebyshev polynomials, \emph{Algebr. Represent. Theory}, {\bf 15}, (2012), no. 3, 527-549.

\bibitem[{FZ02}]{}  { S. Fomin and A. Zelevinsky},
 Cluster algebras I. Foundations, \emph{J. Amer. Math. Soc.}
{\bf 15} (2002), no. 2, 497-529 (electronic) 

\bibitem[FZ03]{} S. Fomin and A. Zelevinsky, Cluster algebras II:
Finite type classification, \emph{Invent. Math.} {\bf 154}, 
(2003), 63-121.

\bibitem[{H88}]{}  { D. Happel}, \emph{Triangulated Categories in the
  Representation Theory of Finite Dimensional Algebras}, London
  Mathematical Society. Lecture Notes Series 119, Cambridge
  University Press, Cambridge, (1988).
  
  \bibitem[KS11] {} B. Keller and S. Scherotzke, Linear recurrence relations for cluster variables of affine quivers \emph{Adv. Math. } {\bf 228} (2011), no. 3, 1842-1862.
  
  \bibitem[{Ma13}]{}  K. E. Magnani, Friezes of type $\mathbb{D}$, submitted for publication.

 \bibitem[{M11}]{}  G. Musiker, A graph theoretic expansion for cluster algebras of classical type, \emph{Ann. Comb.} {\bf 15}, (2011), no. 1, 147-184.
  
\bibitem[{Pr08}]{} J. Propp, The combinatorics of frieze patterns and Markoff numbers, (2008). arXiv:math/0511633v4 [math.CO].

\bibitem[SS] {} D. Simson, A. Skowro\'nski, Elements of the Representation Theory of Associative Algebras, Volume 2: Tubes and Concealed Algebras of Euclidean Type,
in:\emph{London Mathematical Society Student Texts}, vol. 71, \emph{Cambridge University Press}, 2007.
  
\bibitem[{S08}]{}  R. Schiffler, A geometric model for cluster categories of type Dn, \emph{J. Algebraic Combin.} 27 (1) (2008), 1-21.  
\end{thebibliography}
\end{document}